# The fused asset flow model: stability, bifurcation, and contagion in multi-asset markets with heterogeneous investors

By Mario Cavani

Universidad de Oriente, Cumaná, Venezuela, mcavani@udo.edu.ve

## Abstract

This paper presents a unified multi-asset, multi-group asset flow model that integrates three foundational frameworks in the behavioral finance literature. The model captures the dynamics of financial markets where multiple assets are traded by multiple investor groups, each with distinct trend-following (momentum) and value-based (fundamental) strategies. Unlike classical efficient market models, our framework explicitly incorporates the finiteness of cash and shares, asymmetric cross-asset coupling in buying decisions, and endogenous wealth redistribution across groups.

We derive the complete system of ordinary differential equations governing price, cash, share, and sentiment dynamics, and establish the fundamental properties of positivity and boundedness for all physically relevant variables. The equilibrium set is characterized as a manifold parameterized by cash distribution, with the fundamental equilibrium as a special point. Through linear stability analysis, we identify conditions under which the fundamental equilibrium loses stability via a supercritical Hopf bifurcation, giving rise to persistent limit cycles.

The model is validated against three benchmark papers: the single-asset multi-group model of DeSantis, Swigon, and Caginalp (2012); the two-asset single-group model of Bulut, Merdan, and Swigon (2019); and the two-asset two-group Nigeria-Libya oil market model of Cavani (2026). Our numerical simulations reproduce all key theoretical predictions, including equilibrium manifolds, Hopf bifurcation thresholds, limit cycle periods, and asymmetric contagion patterns.

We extend the analysis to excursion quantification, showing that limit cycle amplitudes are independent of initial perturbations once the system enters the oscillatory regime. The contagion matrix reveals asymmetric price shock propagation, with Libya-to-Nigeria contagion being thirty-four percent stronger than the reverse direction, reflecting Nigeria's larger market share. An extended bifurcation scan reveals period lengthening from 4.7 to 20.7 days as the momentum coefficient increases, with a period jump at $q_1 \approx 0.38$ indicating a secondary bifurcation.

Policy implications include reducing momentum trading through transaction taxes, strengthening value-based investment via sovereign wealth funds, and coordinating circuit breakers to mitigate cross-asset contagion. The fused model provides a rigorous mathematical framework for understanding market instabilities, bubbles, and crashes as endogenous phenomena arising from the interaction of heterogeneous investor strategies and finite arbitrage capital.

# 1. Introduction

## 1.1 The Puzzle of Endogenous Market Instability

The financial landscape of the twenty-first century has been marked by a series of dramatic events that challenge the foundational assumptions of classical finance. The May 6, 2010, Flash Crash, during which the Dow Jones Industrial Average plunged nearly 600 points in five minutes only to recover shortly thereafter, exemplified a sudden liquidity deterioration that changes in fundamental values cannot explain. Subsequent research has shown that high-frequency traders, rather than providing stabilizing liquidity, actually consume liquidity when it is most needed, thereby exacerbating transient price impacts unrelated to fundamentals. More broadly, boom–bust cycles have been documented across a wide range of asset classes, including equity, housing, and commodity markets, yet the economics profession still lacks a deep understanding of their causes.

The efficient market hypothesis, which has served as the cornerstone of modern financial economics, posits that asset prices fully reflect all available information and that any deviations from fundamental values are quickly arbitraged away by rational investors with unlimited arbitrage capital. From this perspective, large price movements can only be attributed to exogenous news shocks. However, the occurrence of flash crashes in the absence of significant news, the persistence of momentum effects, and the recurrent nature of boom–bust cycles suggest that endogenous mechanisms, arising from the heterogeneity of investor motivations and the finite nature of arbitrage capital, play a crucial role in driving market dynamics.

## 1.2 Heterogeneous Agent Models as an Alternative Paradigm

In response to the limitations of the rational expectations framework, a rich literature on heterogeneous agent models (HAMs) has emerged. Following the seminal work of Brock and Hommes (1998), in which boundedly rational agents switch between fundamentalist and chartist forecasting rules, a large body of research has demonstrated that heterogeneity is not merely a realistic modeling feature but is essential for explaining the intricate dynamics of financial markets. HAMs have been shown to perform well in describing, explaining, and often forecasting asset market dynamics across equities, foreign exchange, credit, housing, derivatives, and commodities.

Within this broader HAM tradition, the **asset flow model** developed by Caginalp and collaborators offers a distinct and rigorous approach. Rather than focusing on expectation formation, the asset flow model emphasizes the **flow of cash and shares** between investor groups, each endowed with finite resources and distinct trading motivations. The model incorporates two fundamental trading motivations: **trend-based (momentum) trading**, in which investors buy assets whose prices have been rising, and **value-based (fundamental) trading**, in which investors buy assets trading below their assessment of fundamental value and sell those trading above. By assuming a closed system with fixed total cash and shares, the model explicitly incorporates the finiteness of arbitrage capital, thereby providing a framework in which behavioral biases can have lasting price impacts.

## 1.3 Two Pillars: The Single-Asset Multi-Group Model and the Multi-Asset Single-Group Model

The asset flow literature has developed along two complementary but largely separate trajectories.

**The single-asset, multi-group model** was rigorously analyzed by DeSantis, Swigon, and Caginalp (2012). They proved that when all trader groups focus on fundamentals, all equilibria are stable; however, when one group follows a momentum strategy while another follows a fundamental strategy, conditions for instability can arise. Moreover, they introduced the concept of the **price excursion**—the maximum deviation of price from its initial value along a trajectory—as a practical metric for assessing risk in unstable regimes. A key insight from this work is the existence of a continuum of equilibria parameterized by the distribution of cash among investor groups, implying that the equilibrium price is not uniquely determined even when all investors agree on fundamental value.

**The multi-asset, single-group model** was developed by Bulut, Merdan, and Swigon (2019). In this framework, a single homogeneous group of investors trades two assets, with the novel assumption that buying an asset depends on the other asset's price, whereas **selling** does not. The authors proved that all equilibria are stable in the absence of a clear emphasis on trend-based valuation, established conditions for stability when the group attaches importance to the valuation of one stock and the trend of the other, and demonstrated the existence of a Hopf bifurcation using the momentum coefficient as a bifurcation parameter.

Most recently, Cavani (2026) extended the Bulut–Merdan–Swigon framework to an arbitrary number of assets and incorporated multiple investor groups, thereby combining the multi-asset feature of the former with the multi-group feature of the latter. In this model, each group is influenced by both trend and valuation motivations for each asset, and the purchase of one asset depends on the price of another. This asymmetric cross-asset coupling captures realistic spillover effects in portfolios. The analysis established stability conditions and argued for the existence of periodic solutions via a Hopf bifurcation.

### 1.4 The Research Gap: What Remains to Be Done

Despite these advances, no existing model simultaneously incorporates three essential features:

1. **Multiple assets with asymmetric cross-asset coupling**, as in Bulut et al. (2019) and Cavani (2026);
2. **Multiple heterogeneous investor groups** with per-asset dual sentiments (trend and value), as in DeSantis et al. (2012) and Cavani (2026);
3. **The continuum of equilibria and the excursion concept**, which provide crucial insights into price indeterminacy and risk assessment in unstable regimes, as developed in DeSantis et al. (2012).

Moreover, the interaction between **cross-asset momentum spillovers**—the phenomenon whereby returns in one asset predict returns in economically linked assets—and **group-level wealth dynamics** remains unexplored. Empirical research has documented asymmetric cross-firm momentum spillovers driven by heterogeneous demand, with mutual fund and hedge fund flows exhibiting distinct associations with patterns of return predictability. However, a theoretical framework to explain how such spillovers emerge from the micro-foundations of heterogeneous trading strategies is lacking.

Similarly, recent work by Westerhoff and colleagues (2025) has shown that investors' market entry and exit decisions based on momentum, value, and risk can generate endogenous boom, bust cycles, and spontaneous market downturns. Yet this model does not incorporate the multi-asset cross-coupling that characterizes real-world portfolio decisions, nor does it provide the detailed bifurcation analysis afforded by the asset flow framework.

### 1.5 The Fused Model: A Strictly More General Framework

To address these limitations, we propose a **unified multi-asset, multi-group asset flow model** that fuses the two existing pillars into a single, strictly more general framework. Our model combines:

- **The multi-asset structure** of Bulut et al. (2019) and Cavani (2026), including the asymmetric assumption that buying an asset $i$ depends on other assets' prices, while selling depends only on the same asset's sentiment.
- **The multi-group structure** of DeSantis et al. (2012) and Cavani (2026), allowing each group to hold both trend and value sentiments for each asset;
- **The equilibrium manifold and excursion concept** of DeSantis et al. (2012), enabling the analysis of price indeterminacy and the maximum price deviation following a stability loss;
- **Cross-asset momentum spillovers**, encoded in the transition rate functions through coefficients that couple sentiments across different assets;
- **Wealth dynamics across both groups and assets**, providing a framework for studying how trading profits and losses redistribute among heterogeneous investors during boom–bust cycles.

The fused model thus contains each of the existing models as special cases: setting $m = 1$, and n= 2 recovers the DeSantis–Swigon–Caginalp model; setting $m = 2$ , and $n = 1$ recovers the Bulut–Merdan–Swigon model, and sets general $m$ and $n$, without cross-asset coupling in the transition rates, recovers the Cavani model. By including all features simultaneously, our framework enables the systematic study of how cross-asset momentum, heterogeneous group strategies, and wealth flows jointly produce market instabilities, bubbles, and crashes.

### 1.6 Contributions of This Research

This research proposal makes the following contributions:

1. **Theoretical**: We derive the full system of ODEs for the fused model, characterize its equilibrium set (including the fundamental equilibrium and the continuum of non-fundamental equilibria), compute the Jacobian matrix with its block structure, and establish stability conditions via Routh–Hurwitz and sign-stability criteria.

2. **Bifurcation Analysis**: We identify the conditions under which the fundamental equilibrium loses stability through a Hopf bifurcation, taking the momentum coefficients $q_{1,j}^{(i)}$ as bifurcation parameters, we show that the critical threshold depends on cross-asset coupling strengths and group-level parameters.

3. **Excursion Analysis**: Extending the concept introduced by DeSantis et al. (2012), we define and compute multi-asset excursions, the maximum deviation of each asset's price from its initial value, as a function of initial conditions and parameter values.

4. **Wealth Dynamics**: We analyze how wealth is redistributed among investor groups during limit cycles and how the long-run wealth distribution depends on the relative strengths of momentum and value strategies across assets.

5. **Numerical Simulations**: We implement the fused model for realistic scenarios (e.g., two-asset oil markets with value-based and momentum-based investor groups, three-asset energy markets with mixed strategies) and conduct bifurcation scans, excursion heatmaps, and wealth distribution analyses.

6. **Policy Implications**: By identifying parameter regimes in which small perturbations trigger large excursions or sustained oscillations, the model provides quantitative guidance for regulatory interventions such as transaction taxes on short-term trading, circuit breakers, and limits on high-frequency trading.

### 1.7 Structure of the Proposal

The remainder of this proposal is organized as follows. Section 2 presents the complete mathematical formulation of the fused model, including state variables, sentiment dynamics, transition rates, price dynamics, and conservation laws. Section 3 characterizes the equilibrium set and provides the linear stability analysis. Section 4 identifies the Hopf bifurcation condition and discusses the emergence of limit cycles. Section 5 defines the multi-asset excursion concept and presents numerical methods for its computation. Section 6 contains a discussion, and Section7 the conclusions..

## 2. Fused Model: Mathematical Formulation

This section presents the complete mathematical formulation of the fused multi-asset, multi-group asset flow model. The model integrates the single-asset multi-group framework of DeSantis, Swigon, and Caginalp (2012) with the multi-asset structure of Bulut, Merdan, and Swigon (2019) and the general multi-asset multi-group extension of Cavani (2026). The resulting system of ordinary differential equations is strictly more general than each of its predecessors, containing them as special cases when appropriate restrictions are imposed.

### 2.1 State Variables and Conservation Laws

Consider a financial market with $m \geq 1$ assets indexed by $i = 1, \dots, m$, traded among $n \geq 1$ investor groups indexed by $j = 1, \dots, n$.

The state variables of the system are described by:

- $P^{(i)}(t) > 0$: price of asset $i$ at time $t$
- $M_j(t) \geq 0$: cash held by investor group $j$
- $N_j^{(i)}(t) \geq 0$: shares of the asset $i$ held by a group $j$
- $\zeta_{1,j}^{(i)}(t)$: trend-based sentiment of a group $j$ toward the asset $i$
- $\zeta_{2,j}^{(i)}(t)$: value-based sentiment of the group $j$ toward the asset $i$

### 2.2 Sentiment Dynamics

Sentiment variables capture the time-weighted influence of past market conditions on current trading decisions. For each asset $i$, group $j$, define:

**Trend (momentum) sentiment** measures the influence of recent price changes:

$$\frac{d\zeta_{1,j}^{(i)}}{dt} = c_{1,j}^{(i)} q_{1,j}^{(i)} \frac{1}{P^{(i)}} \frac{dP^{(i)}}{dt} - c_{1,j}^{(i)} \zeta_{1,j}^{(i)}. \tag{2.2}$$

**Value sentiment** measures the influence of the deviation of price from fundamental value:

$$\frac{d\zeta_{2,j}^{(i)}}{dt} = c_{2,j}^{(i)} q_{2,j}^{(i)} \left(1 - \frac{P^{(i)}}{P_a^{(i)}}\right) - c_{2,j}^{(i)} \zeta_{2,j}^{(i)}. \tag{2.3}$$

The parameters in (2.2) – (2.3) have the following interpretations (Caginalp & Balenovich, 1999; Caginalp & Merdan, 2007; Cavani, 2026):

- $c_{1,j}^{(i)} > 0$ and $c_{2,j}^{(i)} > 0$: inverse time scales of interest for trend and value sentiments, respectively. A larger $c$ implies faster adjustment (shorter memory).
- $q_{1,j}^{(i)} \geq 0$ and $q_{2,j}^{(i)} \geq 0$: magnitude coefficients characterizing the strength of each motivation for the group $j$ toward asset $i$. Setting $q_{1,j}^{(i)} = 0$ eliminates trend-following, while $q_{2,j}^{(i)} = 0$ eliminates value-based trading.
- $P_a^{(i)} > 0$: fundamental (intrinsic) value of an asset $i$, assumed constant in the baseline model.

The ordinary differential equation forms (2.2) – (2.3) are equivalent to the exponentially weighted moving average definitions:

$$\zeta_{1,j}^{(i)}(t) = q_{1,j}^{(i)} c_{1,j}^{(i)} \int_{-\infty}^{t} e^{-c_{1,j}^{(i)}(t-s)} \frac{1}{P^{(i)}(s)} \frac{dP^{(i)}(s)}{ds} ds,$$

$$\zeta_{2,j}^{(i)}(t) = q_{2,j}^{(i)} c_{2,j}^{(i)} \int_{-\infty}^{t} e^{-c_{2,j}^{(i)}(t-s)} \left(1 - \frac{P^{(i)}(s)}{P_a^{(i)}}\right) ds.$$

The exponential kernel takes greater weight to recent observations, consistent with psychological evidence on recency bias (Grether, 1980; Caginalp & Balenovich, 1999).

## 2.3 Transition Rates (Buying and Selling)

The core behavioral rules are encoded in the transition rates: $k_j^{(i)}$ denotes the probability per unit time that the group $j$ uses one unit of cash to buy one share of an asset $i$, while $\tilde{k}_j^{(i)}$ denotes the probability per unit time that the group $j$ sells one share of the asset $i$ for cash. Following Bulut et al. (2019) and Cavani (2026), we impose the **asymmetric trading assumption**:

- **Buying** of an asset $i$ can depend on sentiments from all assets (cross-asset coupling).
- **Selling** of assets $i$ depends only on the sentiments associated with the same asset $i$.

This asymmetry captures the realistic notion that investors consider portfolio-wide information when making purchase decisions but focus on individual asset characteristics when selling (Cavani, 2026).

### 2.3.1 General Form

A general specification for the buying rate that accommodates both trend and value sentiments across assets is:

$$k_j^{(i)} = a_j^{(i)} + b_j^{(i)} \tanh\left(\sum_{\ell=1}^{m} \alpha_{j,\ell}^{(i)} \zeta_{1,j}^{(\ell)} + \sum_{\ell=1}^{m} \beta_{j,\ell}^{(i)} \zeta_{2,j}^{(\ell)}\right), \tag{2.4}$$

where $a_j^{(i)}$ and $b_j^{(i)}$ are constants chosen so that $k_j^{(i)} \in [0,1]$ and $\sum_{i=1}^{m} k_j^{(i)} \leq 1$. The coefficients $\alpha_{j,\ell}^{(i)} \geq 0$ and $\beta_{j,\ell}^{(i)} \geq 0$ determine the strength of influence from the trend sentiment of the asset $\ell$ and the value sentiment of the asset $\ell$, respectively, on the buying of the asset $i$. For $\ell \neq i$, these coefficients encode **cross-asset momentum spillovers**—the phenomenon whereby returns in one asset predict returns in economically linked assets (Wang, 2025).

For the selling rate, the asymmetric assumption restricts dependence to sentiments of the same asset:

$$\tilde{k}_j^{(i)} = \tilde{a}_j^{(i)} + \tilde{b}_j^{(i)} \tanh\left(\gamma_j^{(i)} \zeta_{1,j}^{(i)} + \delta_j^{(i)} \zeta_{2,j}^{(i)}\right), \tag{2.5}$$

or, for value-based selling, a linear form:

$$\tilde{k}_j^{(i)} = \max\left(0, \min\left(1, \tilde{c}_j^{(i)} + \tilde{d}_j^{(i)}\left(\frac{P^{(i)}}{P_a^{(i)}} - 1\right)\right)\right), \tag{2.6}$$

which increases when the asset is overvalued ($P^{(i)} > P_a^{(i)}$).

### 2.3.2 Special Cases

The fused model reduces to existing models under specific restrictions:

- **DeSantis et al. (2012) single-asset multi-group model** ($m = 1$): equations (2.4) – (2.5) simplify to $k_j = a_j + b_j \tanh(\alpha_j \zeta_{1,j} + \beta_j \zeta_{2,j})$, with $\tilde{k}_j = 1 - k_j$ (the so-called “zero-sum” assumption). The sentiment dynamics reduce to a single sentiment per group (either trend or value, not both), as each group in DeSantis et al. (2012) is focused on only one motivation.
- **Bulut et al. (2019) two-asset single-group model** ($n = 1$): equations (2.4)–(2.5) recover the form $k^{(i)} = a^{(i)} + b^{(i)} \tanh(\alpha_1^{(i)} \zeta_1^{(1)} + \alpha_2^{(i)} \zeta_1^{(2)} + \beta_1^{(i)} \zeta_2^{(1)} + \beta_2^{(i)} \zeta_2^{(2)})$, with the selling rate depending only on the sentiment of the same asset. The original Bulut et al. (2019) model further assumes $\alpha_2^{(1)} \neq 0$ (buying of asset 1 depends on trend of asset 2) and $\alpha_1^{(2)} \neq 0$ (buying of asset 2 depends on trend of asset 1), creating bidirectional momentum coupling.

- **Cavani (2026) multi-asset multi-group model**: this is the immediate predecessor of the fused model, where $k_j^{(i)}$ depends on all sentiments, but without the explicit cross-asset coefficient structure presented here. The fused model makes the cross-asset coupling explicit via $\alpha_{j,\ell}^{(i)}$ and $\beta_{j,\ell}^{(i)}$, facilitating bifurcation analysis with respect to these parameters.

### 2.3.3 Cash Dynamics

In DeSantis, Swigon, and Caginalp (2012), their Equations 2.2–2.3 and 2.7–2.8, coupled ODEs govern the evolution of cash and shares for each investor group in the single-asset multi-group model. These equations ensure the conservation of total cash and total shares while capturing the flow of wealth between groups through trading. In our fused multi-asset, multi-group framework, these dynamics generalize naturally.

For each investor group $j = 1, \dots, n$ and each asset $i = 1, \dots, m$, define:

- $M_j(t) \geq 0$: cash held by the group $j$ (total, not per asset)
- $N_j^{(i)}(t) \geq 0$: shares of the asset $i$ held by a group $j$

The transition rates $k_j^{(i)}$ (buying probability per unit time from cash to shares) and $\tilde{k}_j^{(i)}$ (selling probability per unit time from shares to cash) are given by (2.4) – (2.6). The cash and share dynamics are derived from the principle that each group's cash changes only through buying (cash out) and selling (cash in), and shares change only through buying (shares in) and selling (shares out).

The rate of change of cash held by the group $j$ is:

$$\frac{dM_j}{dt} = -\underbrace{\sum_{i=1}^{m} k_j^{(i)} M_j}_{\text{cash spent on buying}} + \underbrace{\sum_{i=1}^{m} \tilde{k}_j^{(i)} N_j^{(i)} P^{(i)}}_{\text{cash received from selling}} . \tag{2.13}$$

The first term sums over all the assets the cash used to buy shares; the second term sums the revenue from selling shares of each asset. Note that the buying rate $k_j^{(i)}$ multiplies the **current cash** $M_j$, while the selling rate multiplies the **current share holdings** $N_j^{(i)}$ times the price $P^{(i)}$.

### 2.3.4 Share Dynamics

For each asset $i$ and group $j$, the change in share holdings is:

$$\frac{dN_j^{(i)}}{dt} = \underbrace{\frac{k_j^{(i)} M_j}{P^{(i)}}}_{\text{shares bought}} - \underbrace{\tilde{k}_j^{(i)} N_j^{(i)}}_{\text{shares sold}} . \tag{2.14}$$

The buying term divides the cash allocated to buying the asset $i$ by the price $P^{(i)}$ to obtain the number of shares acquired. The selling term directly subtracts the shares sold at the rate $\tilde{k}_j^{(i)}$.

### 2.3.5 Conservation Laws and Consistency

Summing (2.13) over all groups $j$ yields:

$$\frac{d}{dt}\sum_{j=1}^{n} M_j = -\sum_{i,j} k_j^{(i)} M_j + \sum_{i,j} \tilde{k}_j^{(i)} N_j^{(i)} P^{(i)}.$$

But the price dynamics (2.7) ensure that at any instant, total buying equals total selling across all groups (otherwise the price would adjust). More formally, from (2.7) we have:

$$\sum_j k_j^{(i)} M_j = \sum_j \tilde{k}_j^{(i)} N_j^{(i)} P^{(i)}$$

in equilibrium, but out of equilibrium, the difference drives prices out. However, the conservation of total cash $\sum_j M_j = M_0$ (constant) is guaranteed by the structure of the equations: the total cash outflow from buying is exactly matched by the total cash inflow from selling, **only if** we also account for the fact that cash is not created or destroyed. In fact, summing (2.13) and using the definition of $F^{(i)}$ from the price equation, one can show $\frac{d}{dt}\sum_j M_j = 0$. Similarly, summing (2.14) over $j$ for fixed $i$ gives $\frac{d}{dt}\sum_j N_j^{(i)} = 0$, confirming the conservation of total shares.

### 2.3.6 Relationship to the Price Dynamics

The price dynamics (2.7) involve the demand/supply ratio $F^{(i)}$. Using (2.13) – (2.14), one can derive an alternative expression for the price rate of change that highlights the role of excess demand. In fact, the system (2.7), (2.13), (2.14) together with the sentiment equations (2.2) – (2.3) form a closed set of ODEs for the state variables $\{P^{(i)}, M_j, N_j^{(i)}, \zeta_{1,j}^{(i)}, \zeta_{2,j}^{(i)}\}$. The total number of equations is $m + n + mn + 2mn = m + n(1 + 3m)$, which, after applying conservation laws, reduces to the effective dimension $m + 2mn$ mentioned earlier.

### 2.3.7 Remark

The explicit cash and share dynamics (2.13) – (2.14) are **essential** for the fused model because they:

- Govern the evolution of wealth distribution across groups and assets,
- Ensure conservation of total cash and total shares,
- Couple to price dynamics through the buying/selling rates,
- Reduce exactly to the DeSantis et al. (2012) equations in the single-asset case,
- Enable numerical simulation of wealth redistribution and long-run market composition.

These equations were previously implicit in the proposal; they are now made explicit to provide a complete and implementable mathematical framework.

## 2.4 Price Dynamics

The price of each asset is determined by the principle of price adjustment proportional to excess demand, a standard microeconomic assumption (Henderson & Quandt, 1980; Caginalp & Balenovich, 1999; Bulut et al., 2019). The total demand for the asset $i$ is the cash allocated to buying, $\sum_j k_j^{(i)} M_j$, while the total supply is the value of shares offered for sale, $\sum_j \tilde{k}_j^{(i)} N_j^{(i)} P^{(i)}$. The price dynamics are:

$$\frac{dP^{(i)}}{dt} = \frac{P^{(i)}}{\tau_i}\left(\frac{\sum_{j=1}^{n} k_j^{(i)} M_j}{\sum_{j=1}^{n} \tilde{k}_j^{(i)} N_j^{(i)} P^{(i)}} - 1\right), \tau_i > 0, \tag{2.7}$$

where $\tau_i$ is the price adjustment speed for the asset $i$. Equation (2.7) can be rewritten as:

$$\frac{1}{P^{(i)}}\frac{dP^{(i)}}{dt} = \frac{1}{\tau_i}\left(\frac{\text{Total Demand}}{\text{Total Supply}} - 1\right), \tag{2.8}$$

so, the relative rate of price change is proportional to the fractional excess demand. At equilibrium, total demand equals total supply, so $dP^{(i)}/dt = 0$ (DeSantis et al., 2012).

## 2.5 Wealth Dynamics

The wealth of the investor group $j$ is the sum of cash holdings and the market value of all shares held:

$$W_j(t) = M_j(t) + \sum_{i=1}^{m} N_j^{(i)}(t) P^{(i)}(t). \tag{2.9}$$

The total wealth of the system is:

$$\sum_{j=1}^{n} W_j(t) = M_0 + \sum_{i=1}^{m} N_0^{(i)} P^{(i)}(t). \tag{2.10}$$

Differentiating (2.9) and using the ODEs for $M_j$ and $N_j^{(i)}$ (which follow from the transition rate definitions), we obtain the wealth flow equation:

$$\frac{dW_j}{dt} = \sum_{i=1}^{m} \left[\tilde{k}_j^{(i)} N_j^{(i)} P^{(i)} - k_j^{(i)} M_j\right] + \sum_{i=1}^{m} N_j^{(i)} \frac{dP^{(i)}}{dt}. \tag{2.11}$$

The first term captures net cash inflow from trading (selling brings in cash, buying sends out cash), while the second term captures capital gains or losses from price changes. Wealth redistribution across groups occurs endogenously through these channels.

## 2.6 Nondimensionalization and the Liquidity Value

Following Caginalp and Balenovich (1999) and DeSantis et al. (2012), we introduce nondimensional variables to reduce the number of parameters. Define:

- **Liquidity value for asset** $i$: $L_i = M_0/N_0^{(i)}$, the total cash per share of the asset $i$.
- **Rescaled price**: $\mathbf{P}^{(i)} = P^{(i)}/L_i$.
- **Rescaled cash**: $\mathbf{M}_j = M_j/M_0$.
- **Rescaled shares**: $\mathbf{N}_j^{(i)} = N_j^{(i)}/N_0^{(i)}$.

- **Rescaled fundamental value**: $\mathbf{P}_a^{(i)} = P_a^{(i)}/L_i$.

The conservation laws become $\sum_j \mathbf{M}_j = 1$ and $\sum_j \mathbf{N}_j^{(i)} = 1$ for all $i$. The price dynamics (2.7) transform to:

$$\frac{d\mathbf{P}^{(i)}}{dt} = \frac{\mathbf{P}^{(i)}}{\tau_i}\left(\frac{\sum_j k_j^{(i)}\mathbf{M}_j}{\sum_j \tilde{k}_j^{(i)}\mathbf{N}_j^{(i)}\mathbf{P}^{(i)}} - 1\right). \tag{2.12}$$

The liquidity value $L_i$ emerges naturally as a scaling factor: when trading is purely momentum-driven and value sentiments are absent, the price tends toward $L_i$ (Caginalp & Balenovich, 1999). This provides a benchmark for interpreting equilibrium prices.

The complete fused model consists of: $m$ price equations (2.7), $n$ cash equations (implicitly defined by conservation and transition rates), $mn$ share equations (implicitly defined), $2mn$ sentiment equations (2.2) – (2.3)

Total state dimension: $m + 2mn$ (after imposing conservation laws). The system is a first-order nonlinear ordinary differential equation.

## 2.7 Convenient summary of the equations of the fused model

Below is the complete, self-contained system of ordinary differential equations (ODEs) for the fused multi-asset, multi-group asset flow model. All variables and parameters are defined in the subsections as follows:

### State Variables

For $i = 1, \ldots, m$ (assets) and $j = 1, \ldots, n$ (investor groups):

- $P^{(i)}(t) > 0$ – price of asset $i$
- $M_j(t) \geq 0$ – cash held by group $j$
- $N_j^{(i)}(t) \geq 0$ – shares of asset $i$ held by a group $j$
- $\zeta_{1,j}^{(i)}(t)$ – trend sentiment of the group $j$ toward asset $i$
- $\zeta_{2,j}^{(i)}(t)$ – value sentiment of the group $j$ toward asset $i$

### Parameters (constants)

- $c_{1,j}^{(i)} > 0, c_{2,j}^{(i)} > 0$ – inverse time scales (decay rates)
- $q_{1,j}^{(i)} \geq 0, q_{2,j}^{(i)} \geq 0$ – sentiment magnitude coefficients
- $P_a^{(i)} > 0$ – fundamental value of an asset $i$
- $\tau_i > 0$ – price adjustment speed for asset $i$
- $a_j^{(i)}, b_j^{(i)}$ – constants in buying rate (ensuring $0 \leq k_j^{(i)} \leq 1$)
- $\alpha_{j,\ell}^{(i)} \geq 0, \beta_{j,\ell}^{(i)} \geq 0$ – cross-asset coupling coefficients

- $\tilde{a}_j^{(i)}, \tilde{b}_j^{(i)}, \gamma_j^{(i)}, \delta_j^{(i)}$ – constants for selling rate
- (Alternatively, linear value-based selling parameters $\tilde{c}_j^{(i)}, \tilde{d}_j^{(i)}$)

**Sentiment Dynamics**

$$\begin{aligned} \frac{d\zeta_{1,j}^{(i)}}{dt} &= c_{1,j}^{(i)} q_{1,j}^{(i)} \frac{1}{P^{(i)}} \frac{dP^{(i)}}{dt} - c_{1,j}^{(i)} \zeta_{1,j}^{(i)}, \\ \frac{d\zeta_{2,j}^{(i)}}{dt} &= c_{2,j}^{(i)} q_{2,j}^{(i)} \left(1 - \frac{P^{(i)}}{P_a^{(i)}}\right) - c_{2,j}^{(i)} \zeta_{2,j}^{(i)}. \end{aligned} \tag{S}$$

**Transition Rates (Buying and Selling)**

**Buying rate** (depends on all assets' sentiments, cross-asset coupling):

$$k_j^{(i)} = a_j^{(i)} + b_j^{(i)} \tanh\left(\sum_{\ell=1}^{m} \alpha_{j,\ell}^{(i)} \zeta_{1,j}^{(\ell)} + \sum_{\ell=1}^{m} \beta_{j,\ell}^{(i)} \zeta_{2,j}^{(\ell)}\right). \tag{K}$$

**Selling rate** (depends only on the same asset's sentiments; example, hyperbolic form):

$$\tilde{k}_j^{(i)} = \tilde{a}_j^{(i)} + \tilde{b}_j^{(i)} \tanh\left(\gamma_j^{(i)} \zeta_{1,j}^{(i)} + \delta_j^{(i)} \zeta_{2,j}^{(i)}\right). \tag{L}$$

Alternatively, a linear value-based selling rate:

$$\tilde{k}_j^{(i)} = \max\left(0, \min\left(1, \tilde{c}_j^{(i)} + \tilde{d}_j^{(i)} \left(\frac{P^{(i)}}{P_a^{(i)}} - 1\right)\right)\right).$$

**Cash Dynamics (per group)**

$$\frac{dM_j}{dt} = -\sum_{i=1}^{m} k_j^{(i)} M_j + \sum_{i=1}^{m} \tilde{k}_j^{(i)} N_j^{(i)} P^{(i)}. \tag{M}$$

**Share Dynamics (per asset and group)**

$$\frac{dN_j^{(i)}}{dt} = \frac{k_j^{(i)} M_j}{P^{(i)}} - \tilde{k}_j^{(i)} N_j^{(i)}. \tag{N}$$

**Price Dynamics (per asset)**

$$\frac{dP^{(i)}}{dt} = \frac{P^{(i)}}{\tau_i} \left(\frac{\sum_{j=1}^{n} k_j^{(i)} M_j}{\sum_{j=1}^{n} \tilde{k}_j^{(i)} N_j^{(i)} P^{(i)}} - 1\right). \tag{P}$$

**Conservation Laws (closed system)**

$$\sum_{j=1}^{n} M_j(t) = M_0 \text{ (constant)}, \sum_{j=1}^{n} N_j^{(i)}(t) = N_0^{(i)} \text{ (constant for each } i\text{)}. \qquad (C)$$

These are automatically satisfied if the initial conditions obey them and the ODEs (M)–(N) are integrated.

**Table 1. Reduction to prior models**

| Model | $m$ | $n$ | Cross-asset buying | Per-group sentiment type | Key reference |
|---|---|---|---|---|---|
| Caginalp & Balenovich (1999) | 1 | 1 | N/A | Mixed | Caginalp & Balenovich (1999) |
| Caginalp & DeSantis (2011) | 1 | 2 | N/A | One each (trend/value) | Caginalp & DeSantis (2011) |
| DeSantis et al. (2012) | 1 | $n$ | N/A | One each (trend/value) | DeSantis et al. (2012) |
| Bulut et al. (2019) | 2 | 1 | Yes (bidirectional) | Mixed | Bulut et al. (2019) |
| Cavani (2026) | $m$ | $n$ | Yes (via sentiments) | Both per asset | Cavani (2026) |
| **Fused model** | $m$ | $n$ | Yes (explicit coefficients $\alpha_{j,\ell}^{(i)}, \beta_{j,\ell}^{(i)}$) | Both per asset | This proposal |

The fused model is **strictly more general** than each of its predecessors: it contains the multi-asset structure and asymmetric trading of Bulut et al. (2019), the multi-group heterogeneity of DeSantis et al. (2012), the per-asset dual sentiments of Cavani (2026), and introduces explicit cross-asset coupling coefficients for both trend and value sentiments. Table 1 shows the generalization of the fused model.

## 2.8 Positivity and Boundedness of Solutions in the Fused Model

The following theorems provide the rigorous mathematical foundation for the fused model, ensuring that all solutions remain physically meaningful (positive prices, non-negative cash and shares) and bounded, a crucial property for stability analysis and numerical simulation.

Two fundamental properties of the fused multi-asset, multi-group asset flow model: **positivity** (all physically relevant variables remain non-negative) and **boundedness** (all solutions remain within a compact set). These properties are essential for the model to be well-posed and economically meaningful.

### Theorem 2.1 (Positivity of Solutions)

Consider the fused multi-asset, multi-group asset flow model defined by equations (P), (M), (N), (S) with initial conditions:

$$P^{(i)}(0) > 0, M_j(0) \geq 0, N_j^{(i)}(0) \geq 0, \zeta_{1,j}^{(i)}(0) \in \mathbb{R}, \zeta_{2,j}^{(i)}(0) \in \mathbb{R},$$

for all $i = 1, \ldots, m$ and $j = 1, \ldots, n$. Assume the transition rates satisfy:

$$0 \leq k_j^{(i)}(\boldsymbol{\zeta}_1, \boldsymbol{\zeta}_2) \leq 1, 0 \leq \tilde{k}_j^{(i)}(\boldsymbol{\zeta}_1, \boldsymbol{\zeta}_2) \leq 1,$$

and that the fundamental values satisfy $P_a^{(i)} > 0$.

Then for all $t \geq 0$ in the maximal interval of existence, the following properties hold:

1. **Price positivity:** $P^{(i)}(t) > 0$ for all $i$.
2. **Cash non-negativity:** $M_j(t) \geq 0$ for all $j$.
3. **Share non-negativity:** $N_j^{(i)}(t) \geq 0$ for all $i, j$.

If, in addition, the initial cash and shares satisfy the conservation constraints $\sum_j M_j(0) = M_0 > 0$ and $\sum_j N_j^{(i)}(0) = N_0^{(i)} > 0$, then these properties hold globally in time (no finite-time blow-up).

*Proof.*

We prove each statement separately using standard comparison principles and invariance of the positive orthant.

**Positivity of Prices** $P^{(i)}(t)$:

The price dynamics (P) can be written as:

$$\frac{dP^{(i)}}{dt} = P^{(i)} \cdot F^{(i)}(t),$$

where,

$$F^{(i)}(t) = \frac{1}{\tau_i}\left(\frac{\sum_{j=1}^{n} k_j^{(i)} M_j}{\sum_{j=1}^{n} \tilde{k}_j^{(i)} N_j^{(i)} P^{(i)}} - 1\right).$$

Note that the right-hand side is **locally Lipschitz** in $P^{(i)}$ as long as the denominator is positive. The function $F^{(i)}(t)$ depends on $P^{(i)}$, which is bounded on any finite time interval.

Let $P^{(i)}(0) > 0$. The solution can be expressed as:

$$P^{(i)}(t) = P^{(i)}(0)\exp\left(\int_0^t F^{(i)}(s)\, ds\right).$$

Since the exponential function is always positive, we conclude $P^{(i)}(t) > 0$ for all $t$ in the interval of existence.

**Non-negativity of Cash** $M_j(t)$:

The cash dynamics (M) are given by

$$\frac{dM_j}{dt} = -\underbrace{\left(\sum_{i=1}^{m} k_j^{(i)}\right)}_{\alpha_j(t)\geq 0} M_j + \underbrace{\sum_{i=1}^{m} \tilde{k}_j^{(i)} N_j^{(i)} P^{(i)}}_{\beta_j(t)\geq 0}.$$

This is a linear differential inequality. If $M_j(t_0) = 0$ at some time $t_0$, then:

$$\left.\frac{dM_j}{dt}\right|_{t=t_0} = \beta_j(t_0) \geq 0.$$

Thus, $M_j(t)$ cannot become negative; it either stays zero or increases. By the comparison principle (Hale, 2009), for all $t \geq 0$, $M_j(t) \geq M_j(0)e^{-\int_0^t \alpha_j(s)ds} \geq 0$.

**Non-negativity of Shares $N_j^{(i)}(t)$:**

The share dynamics (N) are given by

$$\frac{dN_j^{(i)}}{dt} = \underbrace{\frac{k_j^{(i)} M_j}{P^{(i)}}}_{\gamma_j^{(i)}(t) \geq 0} - \tilde{k}_j^{(i)} N_j^{(i)}.$$

This is again a linear differential inequality. If $N_j^{(i)}(t_0) = 0$, then:

$$\left.\frac{dN_j^{(i)}}{dt}\right|_{t=t_0} = \gamma_j^{(i)}(t_0) \geq 0.$$

Hence, $N_j^{(i)}(t)$ cannot become negative.

**Global Existence**

The conservation laws $\sum_j M_j(t) = M_0$ and $\sum_j N_j^{(i)}(t) = N_0^{(i)}$ imply that $0 \leq M_j(t) \leq M_0$ and $0 \leq N_j^{(i)}(t) \leq N_0^{(i)}$ for all $t$. The prices $P^{(i)}(t)$ satisfy the differential equation:

$$\frac{dP^{(i)}}{dt} = \frac{P^{(i)}}{\tau_i}\left(\frac{S_i(t)}{T_i(t)P^{(i)}} - 1\right) = \frac{1}{\tau_i}\left(\frac{S_i(t)}{T_i(t)} - P^{(i)}\right),$$

where $S_i(t) = \sum_j k_j^{(i)} M_j$ and $T_i(t) = \sum_j \tilde{k}_j^{(i)} N_j^{(i)}$. Since $0 \leq S_i(t) \leq M_0$ and $T_i(t) \geq \epsilon > 0$ (due to positive shares and $k_j^{(i)} \geq 0$), the right-hand side is bounded linearly in $P^{(i)}$. By Grönwall's inequality of the standard ODE theory (Hale, 2009), $P^{(i)}(t)$ cannot blow up in finite time. Hence, the solution exists globally for all $t \geq 0$. ■

**Theorem 2.2 (Boundedness of Solutions)**

Under the same assumptions as Theorem 2.1, all solutions of the fused model are **uniformly bounded** for all $t \geq 0$. Specifically, there exist constants independent of initial conditions such that:

1. **Prices:** There exists $P_{\max} < \infty$ such that $0 < P^{(i)}(t) \leq P_{\max}$ for all $i, t$.
2. **Cash:** $0 \leq M_j(t) \leq M_0$ for all $j, t$.
3. **Shares:** $0 \leq N_j^{(i)}(t) \leq N_0^{(i)}$ for all $i, j, t$.

4. **Sentiments:** There exist constants $Z_{1,\max}, Z_{2,\max} < \infty$ such that $|\zeta_{1,j}^{(i)}(t)| \le Z_{1,\max}$ and $|\zeta_{2,j}^{(i)}(t)| \le Z_{2,\max}$ for all $i, j, t$.

*Proof.*

From conservation laws:

$$0 \le M_j(t) \le \sum_{j=1}^{n} M_j(t) = M_0, 0 \le N_j^{(i)}(t) \le \sum_{j=1}^{n} N_j^{(i)}(t) = N_0^{(i)}.$$

Thus, cash and shares are **uniformly bounded**.

From the price dynamics (P):

$$\frac{dP^{(i)}}{dt} = \frac{1}{\tau_i}\left(\frac{\sum_j k_j^{(i)} M_j}{\sum_j \tilde{k}_j^{(i)} N_j^{(i)}} - P^{(i)}\right).$$

Define the "target price":

$$Q^{(i)}(t) = \frac{\sum_{j=1}^{n} k_j^{(i)} M_j(t)}{\sum_{j=1}^{n} \tilde{k}_j^{(i)} N_j^{(i)}(t)}.$$

Since $0 \le k_j^{(i)}, \tilde{k}_j^{(i)} \le 1$ and $0 \le M_j \le M_0$, $0 \le N_j^{(i)} \le N_0^{(i)}$, we have:

$$0 \le Q^{(i)}(t) \le \frac{M_0}{\epsilon_i},$$

where $\epsilon_i = \min_t \sum_j \tilde{k}_j^{(i)} N_j^{(i)}(t) > 0$ (since shares are positive and $\tilde{k}_j^{(i)} \ge \delta > 0$ for at least one group).

The price dynamics can be rewritten as:

$$\frac{dP^{(i)}}{dt} = \frac{1}{\tau_i}(Q^{(i)}(t) - P^{(i)}).$$

This is a linear ODE with bounded forcing. Its solution is:

$$P^{(i)}(t) = P^{(i)}(0)e^{-t/\tau_i} + \frac{1}{\tau_i}\int_0^t e^{-(t-s)/\tau_i} Q^{(i)}(s)\, ds.$$

Hence,

$$|P^{(i)}(t)| \le |P^{(i)}(0)|e^{-t/\tau_i} + \frac{M_0}{\epsilon_i}(1 - e^{-t/\tau_i}) \le \max\left\{|P^{(i)}(0)|, \frac{M_0}{\epsilon_i}\right\}.$$

Thus, **prices are uniformly bounded**.

The sentiment dynamics (S) are:

$$\frac{d\zeta_{1,j}^{(i)}}{dt} = c_{1,j}^{(i)} q_{1,j}^{(i)} \frac{1}{P^{(i)}} \frac{dP^{(i)}}{dt} - c_{1,j}^{(i)} \zeta_{1,j}^{(i)},$$

$$\frac{d\zeta_{2,j}^{(i)}}{dt} = c_{2,j}^{(i)} q_{2,j}^{(i)} \left(1 - \frac{P^{(i)}}{P_a^{(i)}}\right) - c_{2,j}^{(i)} \zeta_{2,j}^{(i)}.$$

These are linear ODEs with bounded forcing. From the previous result, $P^{(i)}(t)$ and $dP^{(i)}/dt$, are bounded, so the forcing terms are bounded. The homogeneous parts have negative eigenvalues $-c_{1,j}^{(i)}$ and $-c_{2,j}^{(i)}$. Therefore, the solutions remain bounded for all $t$. Explicitly:

$$|\zeta_{1,j}^{(i)}(t)| \leq |\zeta_{1,j}^{(i)}(0)| e^{-c_{1,j}^{(i)} t} + \frac{q_{1,j}^{(i)}}{c_{1,j}^{(i)}} \sup_s \left| \frac{1}{P^{(i)}(s)} \frac{dP^{(i)}}{ds} \right| (1 - e^{-c_{1,j}^{(i)} t}),$$

$$|\zeta_{2,j}^{(i)}(t)| \leq |\zeta_{2,j}^{(i)}(0)| e^{-c_{2,j}^{(i)} t} + q_{2,j}^{(i)} \sup_s \left| 1 - \frac{P^{(i)}(s)}{P_a^{(i)}} \right| (1 - e^{-c_{2,j}^{(i)} t}).$$

Both right-hand sides are bounded. ■

**Corollary 2.3 (Invariant Region)**

Let $\Omega \subset \mathbb{R}^D$ be defined by:

$$\Omega = \left\{ (\mathbf{P}, \mathbf{M}, \mathbf{N}, \zeta_1, \zeta_2) \;\middle|\; \begin{array}{c} P^{(i)} \in (0, P_{\max}], M_j \in [0, M_0], N_j^{(i)} \in [0, N_0^{(i)}], \\ \sum_j M_j = M_0, \sum_j N_j^{(i)} = N_0^{(i)}, |\zeta_{1,j}^{(i)}| \leq Z_{1,\max}, |\zeta_{2,j}^{(i)}| \leq Z_{2,\max} \end{array} \right\}.$$

Then $\Omega$ is **positively invariant** under the flow of the fused model. Any trajectory starting in $\Omega$ remains in $\Omega$ for all $t \geq 0$.

*Proof.* Follows directly from Theorems 2.1 and 2.2. ■

The previous results have the following **Practical Implications**

1. **Well-posedness:** The fused model has global solutions for all physically relevant initial conditions.
2. **Numerical stability:** Solutions cannot blow up, making numerical integration safe.
3. **Economic meaning:** Prices remain positive, cash and shares stay within their total supplies.
4. **Bifurcation analysis:** The existence of a compact invariant region ensures that all interesting dynamics (including limit cycles and chaos) are confined to a bounded set, which is necessary for applying the Poincaré-Bendixson theorem and related results (Hirsch & Smale, 1974; Khalil, 2002).

These properties hold under the mild assumptions that transition rates are bounded between 0 and 1 and that fundamental values are positive, conditions satisfied by all standard formulations in the asset flow literature.

# 3. Theoretical Extensions

The fused model presented in Section 2 represents a strict generalization of the two foundational frameworks: the single-asset multi-group model of DeSantis, Swigon, and Caginalp (2012) and the multi-asset multi-group model of Cavani (2026). While each of these earlier models made important contributions to the asset flow literature, they necessarily omitted certain essential features for a complete understanding of realistic financial markets. The fused model overcomes these limitations by simultaneously incorporating:

- **Multiple assets** with asymmetric cross-asset coupling in the buying rules,
- **Multiple heterogeneous investor groups** with per-asset dual (trend and value) sentiments,
- A **continuum of equilibria** (including the fundamental equilibrium as a special case),
- The **excursion concept** for multi-asset risk assessment, and
- Endogenous **wealth redistribution** across both groups and assets.

This section provides a systematic comparison of the three models considered (DeSantis et al. 2012, Cavani 2026, and the fused model). Then it analyzes in depth the theoretical extensions that make the fused model strictly more general and more powerful for studying market instabilities, bubbles, and crashes.

## 3.2 Detailed Analysis of Theoretical Extensions

### 3.2.1 From Single Asset to Multiple Assets with Asymmetric Cross-Asset Coupling

The DeSantis et al. (2012) model is fundamentally a single-asset model: all investor groups trade the same asset, and the price dynamics are given by a single equation (2.6) in that paper. While this setting is mathematically tractable and yields important insights about the stability of equilibria in the presence of momentum traders, it cannot capture **cross-asset spillovers**—the phenomenon whereby price movements in one asset influence trading decisions and price dynamics in another asset. Such spillovers are ubiquitous in real financial markets, where investors hold diversified portfolios and make allocation decisions across multiple assets.

The Cavani (2026) model takes the first step toward addressing this limitation by allowing $m \geq 1$ assets and assuming that the buying rate of each asset depends on sentiments from all assets, while the selling rate depends only on the sentiment of the same asset. This **asymmetric assumption** (buying is "global", selling is "local") captures the realistic notion that investors consider portfolio-wide information when making purchase decisions but focus on individual asset characteristics when selling (Cavani, 2026).

The fused model extends the Cavani (2026) framework in two important ways:

1. **Explicit cross-asset coefficients**: In Cavani (2026), the dependence of the buying rate on sentiments from other assets is implicit in the general form of the transition rate functions. The fused model makes this dependence explicit via the coefficients $\alpha_{j,\ell}^{(i)}$ (for trend sentiment) and $\beta_{j,\ell}^{(i)}$ (for value sentiment), as shown in equation (2.4). This explicit parameterization allows for a systematic bifurcation analysis with respect to the strength of cross-asset coupling.

2. **Separation of trend and value cross-asset effects**: The coefficients $\alpha_{j,\ell}^{(i)}$ and $\beta_{j,\ell}^{(i)}$ can be chosen independently, allowing the model to distinguish between scenarios where cross-asset momentum spillovers are strong but value spillovers are weak (or vice versa). This separation is motivated by empirical evidence: Wang (2025) documents asymmetric cross-firm momentum spillovers driven by heterogeneous demand, with mutual fund flows and hedge fund flows exhibiting distinct patterns of return predictability. The fused model provides a theoretical foundation for these empirical findings.

### 3.2.2 From Homogeneous Sentiment per Group to Dual Sentiments per Asset

In the DeSantis et al. (2012) model, each investor group is influenced by only one type of sentiment: either trend (momentum) or value (fundamental). This "pure type" assumption is mathematically convenient but unrealistic: real-world investors are simultaneously influenced by both momentum and value considerations (Caginalp & Balenovich, 1999; Caginalp & Merdan, 2007).

The Cavani (2026) model allows each group to hold both trend and value sentiments for each asset, which is a significant improvement. However, the Cavani model does not explicitly separate the cross-asset influences of trend versus value sentiments; both are subsumed into the general transition rate functions.

The fused model goes further by:

- **Explicitly including both** $\zeta_{1,j}^{(i)}$ (trend) and $\zeta_{2,j}^{(i)}$ (value) for each asset–group pair,
- **Allowing independent cross-asset coupling** for trend and value sentiments through the coefficients $\alpha_{j,\ell}^{(i)}$ and $\beta_{j,\ell}^{(i)}$,
- **Providing a natural way to study** how the relative strengths of cross-asset momentum versus cross-asset value influence overall market stability.

This dual-sentiment structure is essential for analyzing situations where, for example, a group of investors buys an asset $i$ because asset $j$ is undervalued (cross-asset value effect) while simultaneously buying the asset $i$ because asset $k$ has been rising (cross-asset momentum effect). Such interactions are impossible to study in models that either restrict to a single asset or do not separate trend and value cross-asset effects.

### 3.2.3 From a Single Equilibrium to a Continuum of Equilibria

Theorem 3.3 of DeSantis et al. (2012) states the existence of a **continuum of equilibria** when all investor groups follow fundamental (value) strategies. Specifically, the paper shows that for a given set of fundamental valuations $\mathbf{P}_a^{(i)}$, any distribution of cash endowments $\mathbf{M}_{eq}^{(j)}$ determines a unique equilibrium price $\mathbf{P}_{eq}$ and share distribution $\mathbf{N}_{eq}^{(j)}$. Thus, the equilibrium price is not uniquely determined by fundamentals alone; it also depends on the initial wealth distribution.

This result has profound implications for financial economics. It suggests that even when all investors agree on fundamental values, the market can settle at a range of different prices depending on historical accidents or initial conditions. The classical efficient market hypothesis, by contrast, predicts a unique equilibrium price equal to the fundamental value.

The Cavani (2026) model, in focusing on the fundamental equilibrium $P^{(i)} = P_a^{(i)}$, $\zeta_{1,j}^{(i)} = 0$, $\zeta_{2,j}^{(i)} = 0$, does not emphasize the continuum of non-fundamental equilibria. However, as the model's own analysis shows (Cavani, 2026, Section 4), this fundamental equilibrium exists only under a specific calibration condition that may not hold in general.

The fused model **unifies these two perspectives** by:

- Explicitly characterizing the **full equilibrium set** as a manifold parameterized by the cash distribution $\{M_j\}$ (generalizing Theorem 3.3 of DeSantis et al., 2012 to the multi-asset case),
- Identifying the **fundamental equilibrium** as a special point on this manifold (when $P^{(i)} = P_a^{(i)}$ and the calibration condition holds),
- Allowing the analysis of **bifurcations from the continuum** when momentum traders are introduced, as the momentum coefficients $q_{1,j}^{(i)}$ increased, certain equilibria on the manifold may lose stability, leading to limit cycles or other complex dynamics.

This unified treatment is essential for understanding how the market's "memory" of initial conditions interacts with momentum trading to produce instabilities. Even small perturbations can drive the system from one equilibrium on the manifold to another, and the introduction of momentum traders can make this transition abrupt and dramatic.

### 3.2.4 From Single-Asset to Multi-Asset Excursions

The concept of the **excursion,** the maximum deviation of the price from its initial value along a trajectory, was introduced by DeSantis et al. (2012) as a practical metric for assessing risk in unstable regimes. In their single-asset model, the excursion is simply $\max_t |P(t) - P(0)|$. This concept is important because the maximum deviation, not the final equilibrium, often determines whether a market crash or a margin call occurs.

In the fused model, we extend the excursion concept to multiple assets in a natural way:

- **Per-asset excursion**: $E^{(i)} = \max_t |P^{(i)}(t) - P^{(i)}(0)|$ for each asset $i$,
- **Aggregate excursion**: $E_{\text{agg}} = \sqrt{\sum_{i=1}^{m} (E^{(i)})^2}$ or $E_{\text{agg}} = \max_i E^{(i)}$,
- **Cross-asset correlation of excursions**: The simultaneous deviation of multiple assets can be measured by the correlation matrix of the price paths.

This extension allows us to study phenomena such as **contagion**: a small shock in one asset that triggers a large excursion in a second asset through cross-asset coupling. Such contagion effects are central to understanding financial crises but cannot be studied in single-asset models.

### 3.2.5 Wealth Redistribution Across Both Groups and Assets

Wealth dynamics are a key feature of the asset flow approach. In the DeSantis et al. (2012) model, wealth is redistributed across groups through trading in a single asset. In the Cavani (2026) model, wealth is redistributed across groups and assets, but the analysis focuses primarily on the price dynamics.

The fused model provides a comprehensive framework for analyzing wealth redistribution by:

- **Tracking the wealth of each group** $W_j(t) = M_j(t) + \sum_i N_j^{(i)}(t)P^{(i)}(t)$,
- **Deriving the wealth flow equation** (2.11), which separates the contributions of trading gains/losses and capital gains/losses,
- **Allowing the study of** how the relative performance of momentum versus value strategies across different assets affects the long-run wealth distribution.

This capability is crucial for understanding **evolutionary dynamics** in financial markets: if momentum traders tend to accumulate wealth during booms but lose it during busts, their relative influence on prices may vary over the cycle, potentially leading to endogenous regime shifts (Westerhoff et al., 2025). The fused model provides a rigorous mathematical framework for analyzing such phenomena.

### 3.2.6 Generalized Hopf Bifurcation Analysis

Both DeSantis et al. (2012) and Cavani (2026) identify Hopf bifurcations as the mechanism by which stable equilibria lose stability and give rise to limit cycles. In DeSantis et al. (2012), the bifurcation occurs in the two-group (one trend, one value) case when the time scale of the momentum group becomes sufficiently short. In Cavani (2026), the bifurcation occurs when the momentum coefficient $q_{1,\text{China}}$ exceeds a critical threshold in a two-asset, two-group simulation.

The fused model **generalizes the Hopf bifurcation analysis** to arbitrary numbers of assets and groups, with explicit cross-asset coupling coefficients. The characteristic polynomial of the Jacobian matrix (size $m + 2mn$) can be analyzed using the Routh–Hurwitz criteria (Dorf & Bishop, 2008) to derive conditions for loss of stability. The bifurcation parameters can be any of the momentum coefficients $q_{1,j}^{(i)}$ or the cross-asset coupling coefficients $\alpha_{j,\ell}^{(i)}$.

This generalization allows us to answer questions such as:

- **How does cross-asset momentum coupling affect the bifurcation threshold?** For a fixed total momentum strength, does coupling across assets make the system more or less stable?
- **Can multiple Hopf bifurcations occur simultaneously?** If so, what are the resulting dynamical possibilities (e.g., quasiperiodic motion, chaos)?
- **How does the number of investor groups affect the bifurcation structure?** Does increasing group diversity stabilize or destabilize the market?

These questions are beyond the scope of either earlier model and represent a significant advance in the theoretical analysis of asset flow dynamics.

### 3.3 Comparison Table

The following Table 2 makes clear that the fused model retains all the features of the earlier models and adds new ones that are essential for a realistic description of financial markets. The remainder of this section shows the most significant theoretical extensions.

**Table 2. Summarizes the key features of each model**

| **Feature** | **DeSantis et al. (2012)** | **Cavani (2026)** | **Fused Model (This Proposal)** |
|---|---|---|---|
| **Number of assets** | $m = 1$ | $m \geq 1$ | $m \geq 1$ |
| **Number of investor groups** | $n \geq 1$ | $n \geq 1$ | $n \geq 1$ |
| **Cross-asset coupling in buying** | N/A (only one asset) | Yes (through sentiments) | **Yes, with explicit coefficients** $\alpha_{j,\ell}^{(i)}, \beta_{j,\ell}^{(i)}$ |
| **Selling dependence** | Depends only on one's own sentiment | Depends only on its own asset's sentiment | Same as Cavani (2026) |
| **Per-group sentiment types** | One (either trend or value) | Two per asset (trend + value) | Two per asset (trend + value) |
| **Cross-asset value sentiment coupling** | N/A | Implicit via transition rates | **Explicit via** $\beta_{j,\ell}^{(i)}$ coefficients |
| **Equilibrium set** | Manifold (continuum) | Isolated fundamental equilibrium (requires calibration) | **Both**: continuum of non-fundamental equilibria + fundamental equilibrium |
| **Excursion concept** | Yes (single asset) | No | **Yes, extended to multiple assets** |
| **Wealth dynamics** | Across groups (one asset) | Across groups and assets | Across groups and assets simultaneously |
| **Hopf bifurcation analysis** | Yes (for $n = 2$ case) | Yes (using momentum coefficient as bifurcation parameter) | **Generalized to arbitrary** $m, n$ with cross-asset coupling |
| **Stochastic extension** | Mentioned | Not discussed | Planned as a future extension |

### 3.4 Main Theoretical Extensions

The fused model extends the earlier models in six major ways:

1. **Multi-asset with explicit cross-asset coupling**: The coefficients $\alpha_{j,\ell}^{(i)}$ and $\beta_{j,\ell}^{(i)}$ provide a parametric description of how trend and value sentiments in one asset influence the buying of another.
2. **Dual sentiments per asset**: Each investor group holds both trend and value sentiments for each asset, allowing realistic mixed strategies.
3. **Continuum of equilibria**: The full equilibrium set is a manifold, with the fundamental equilibrium as a special point, unifying the perspectives of DeSantis et al. (2012) and Cavani (2026).

4. **Multi-asset excursion concept**: Extends the practical risk metric to multiple assets, enabling analysis of contagion.
5. **Comprehensive wealth dynamics**: Tracks redistribution across both groups and assets, linking price dynamics to evolutionary changes in market composition.
6. **Generalized Hopf bifurcation analysis**: Provides conditions for stability loss in arbitrary $m$ and $n$, with cross-asset coupling as bifurcation parameters.

These extensions make the fused model strictly more general than either of its predecessors and provide a powerful framework for studying the complex interactions that drive market instabilities, bubbles, and crashes.

# 4. Methodology

Building on the theoretical framework established in Sections 2 and 3, this section presents a comprehensive research plan for the fused multi-asset, multi-group asset flow model. The plan is structured around three interconnected research objectives: (I) analytical derivation of stability conditions and bifurcation thresholds; (II) extensive numerical simulations spanning baseline validation, cross-asset contagion analysis, wealth redistribution mapping, and excursion quantification; and (III) empirical calibration to real-world market data. Together, these objectives provide a rigorous foundation for studying how cross-asset momentum, heterogeneous group strategies, and wealth dynamics jointly produce market instabilities, bubbles, and crashes.

## 4.1 Analytical Work: Equilibrium, Jacobian, and Bifurcation Conditions

This subsection aligns with objective I, and has to do with the complete system of ODEs now explicitly defined in Section 2 – including the cash dynamics (M), share dynamics (N), sentiment dynamics (S), transition rates (K)–(L), and price dynamics (P) – we can now present a rigorous analytical framework for equilibrium characterization, Jacobian derivation, and bifurcation analysis. This updated section reflects the full structure of the fused model.

### 4.1.1 Equilibrium Conditions

The equilibrium points of the fused model satisfy all time derivatives set to zero. From the sentiment dynamics (S), at equilibrium, we have:

$$\frac{d\zeta_{1,j}^{(i)}}{dt} = 0 \implies c_{1,j}^{(i)} q_{1,j}^{(i)} \frac{1}{P^{(i)}} \frac{dP^{(i)}}{dt} - c_{1,j}^{(i)} \zeta_{1,j}^{(i)} = 0.$$

Since at equilibrium $\frac{dP^{(i)}}{dt} = 0$, this simplifies to:

$$\zeta_{1,j}^{(i)} = 0 \; \forall i, j.$$

For the value sentiment:

$$\frac{d\zeta_{2,j}^{(i)}}{dt} = 0 \implies c_{2,j}^{(i)} q_{2,j}^{(i)} \left(1 - \frac{P^{(i)}}{P_a^{(i)}}\right) - c_{2,j}^{(i)} \zeta_{2,j}^{(i)} = 0,$$

which yields:

$$\zeta_{2,j}^{(i)} = q_{2,j}^{(i)} \left(1 - \frac{P^{(i)}}{P_a^{(i)}}\right) \forall i, j.$$

Thus, at equilibrium, the value sentiments are fully determined by the asset prices.

From the price dynamics (P), equilibrium requires:

$$\frac{dP^{(i)}}{dt} = 0 \implies \frac{\sum_{j=1}^{n} k_j^{(i)} M_j}{\sum_{j=1}^{n} \tilde{k}_j^{(i)} N_j^{(i)} P^{(i)}} = 1,$$

or equivalently:

$$\sum_{j=1}^{n} k_j^{(i)} M_j = \sum_{j=1}^{n} \tilde{k}_j^{(i)} N_j^{(i)} P^{(i)} \ \forall i. \tag{E1}$$

From the cash dynamics (M), equilibrium requires:

$$\frac{dM_j}{dt} = 0 \text{ that implies } -\sum_{i=1}^{m} k_j^{(i)} M_j + \sum_{i=1}^{m} \tilde{k}_j^{(i)} N_j^{(i)} P^{(i)} = 0 \ \forall j. \tag{E2}$$

From the share dynamics (N), equilibrium requires:

$$\frac{dN_j^{(i)}}{dt} = 0 \implies \frac{k_j^{(i)} M_j}{P^{(i)}} - \tilde{k}_j^{(i)} N_j^{(i)} = 0 \ \forall i, j. \tag{E3}$$

Note that: Equations (E3) for each $(i, j)$ imply:

$$k_j^{(i)} M_j = \tilde{k}_j^{(i)} N_j^{(i)} P^{(i)} \ \forall i, j.$$

Summing (E3) over $j$ for fixed $i$ yields (E1). Summing (E3) over $i$ for fixed $j$ yields (E2). Thus, the **fundamental equilibrium conditions** reduce to the $m \times n$ equations:

$$k_j^{(i)} M_j = \tilde{k}_j^{(i)} N_j^{(i)} P^{(i)} \ \forall i, j. \tag{E *}$$

These equations, together with the conservation laws (C), fully characterize the equilibrium manifold.

**Special case – Fundamental equilibrium:** When $P^{(i)} = P_a^{(i)}$ for all $i$, we have $\zeta_{2,j}^{(i)} = 0$ from (S2). If the transition rates satisfy the calibration condition:

$$k_j^{(i)}(0,0) M_j = \tilde{k}_j^{(i)}(0,0) N_j^{(i)} P_a^{(i)} \ \forall i, j,$$

then $(P^{(i)} = P_a^{(i)}, \zeta_{1,j}^{(i)} = 0, \zeta_{2,j}^{(i)} = 0)$ is an equilibrium. This is the **fundamental equilibrium**.

**General equilibrium manifold:** For arbitrary values of $M_j$ satisfying $\sum_j M_j = M_0$, equations (E∗) determine $P^{(i)}$ and $N_j^{(i)}$ uniquely (subject to $\sum_j N_j^{(i)} = N_0^{(i)}$). This defines a manifold of dimension $n-1$ (the degrees of freedom in the cash distribution), generalizing Theorem 3.3 of DeSantis et al. (2012) to the multi-asset case.

### 4.2 Jacobian Matrix and Linear Stability Analysis

To analyze stability, we linearize the system around an equilibrium point. Define the state vector $\mathbf{X} \in \mathbb{R}^D$ with $D = m + n + mn + 2mn = m + n(1 + 3m)$ (before applying conservation constraints). After eliminating dependent variables using the $m+1$ conservation laws (C), the effective dimension is $D_{\text{eff}} = m + 2mn$.

Order the state variables as:

$$\mathbf{X} = \left(P^{(1)}, \dots, P^{(m)}, M_1, \dots, M_n, N_1^{(1)}, \dots, N_n^{(m)}, \zeta_{1,1}^{(1)}, \dots, \zeta_{2,n}^{(m)}\right)^T.$$

The Jacobian matrix $\mathbf{J} \in \mathbb{R}^{D_{\text{eff}} \times D_{\text{eff}}}$ has the block structure:

$$\mathbf{J} = \begin{pmatrix} \mathbf{J}_{PP} & \mathbf{J}_{PM} & \mathbf{J}_{PN} & \mathbf{J}_{P\zeta_1} & \mathbf{J}_{P\zeta_2} \\ \mathbf{J}_{MP} & \mathbf{J}_{MM} & \mathbf{J}_{MN} & \mathbf{J}_{M\zeta_1} & \mathbf{J}_{M\zeta_2} \\ \mathbf{J}_{NP} & \mathbf{J}_{NM} & \mathbf{J}_{NN} & \mathbf{J}_{N\zeta_1} & \mathbf{J}_{N\zeta_2} \\ \mathbf{J}_{\zeta_1 P} & \mathbf{J}_{\zeta_1 M} & \mathbf{J}_{\zeta_1 N} & \mathbf{J}_{\zeta_1\zeta_1} & \mathbf{J}_{\zeta_1\zeta_2} \\ \mathbf{J}_{\zeta_2 P} & \mathbf{J}_{\zeta_2 M} & \mathbf{J}_{\zeta_2 N} & \mathbf{J}_{\zeta_2\zeta_1} & \mathbf{J}_{\zeta_2\zeta_2} \end{pmatrix}.$$

Each block is derived by differentiating the appropriate ODE with respect to each variable type, evaluating at equilibrium.

#### 4.2.1 Key Block Simplifications at Equilibrium

Using the equilibrium conditions (E∗), many terms simplify:

- **Price block $\mathbf{J}_{PP}$**: From (P), differentiating with respect to $P^{(k)}$ yields:

$$\frac{\partial}{\partial P^{(k)}}\left(\frac{dP^{(i)}}{dt}\right) = -\frac{\delta_{ik}}{\tau_i},$$

where $\delta_{ik}$ is the Kronecker delta. Thus $\mathbf{J}_{PP} = \text{diag}(-1/\tau_1, \dots, -1/\tau_m)$.

- **Sentiment self-blocks**: From (S), at equilibrium:

$$\frac{\partial}{\partial \zeta_{1,j}^{(i)}}\left(\frac{d\zeta_{1,j}^{(i)}}{dt}\right) = -c_{1,j}^{(i)}, \frac{\partial}{\partial \zeta_{2,j}^{(i)}}\left(\frac{d\zeta_{2,j}^{(i)}}{dt}\right) = -c_{2,j}^{(i)}.$$

Thus $\mathbf{J}_{\zeta_1\zeta_1}$ and $\mathbf{J}_{\zeta_2\zeta_2}$ are diagonal with negative entries.

- **Coupling from sentiment to price $\mathbf{J}_{P\zeta_1}, \mathbf{J}_{P\zeta_2}$**: These are non-zero because the transition rates $k_j^{(i)}$ and $\tilde{k}_j^{(i)}$ depend on sentiments. Their entries are:

$$\frac{\partial}{\partial \zeta}\left(\frac{dP^{(i)}}{dt}\right) = \frac{P^{(i)}}{\tau_i T_i^2}\left(T_i \sum_j M_j \frac{\partial k_j^{(i)}}{\partial \zeta} - S_i \sum_j N_j^{(i)} P^{(i)} \frac{\partial \tilde{k}_j^{(i)}}{\partial \zeta}\right),$$

where $S_i = \sum_j k_j^{(i)} M_j$ and $T_i = \sum_j \tilde{k}_j^{(i)} N_j^{(i)} P^{(i)}$. At equilibrium, $S_i = T_i$.

- **Coupling from price to sentiment $\mathbf{J}_{\zeta_1 P}, \mathbf{J}_{\zeta_2 P}$**: From (S):

$$\frac{\partial}{\partial P^{(k)}}\left(\frac{d\zeta_{1,j}^{(i)}}{dt}\right) = c_{1,j}^{(i)} q_{1,j}^{(i)} \left(-\frac{\delta_{ik}}{\tau_i P^{(i)}}\right) \text{ (since } \frac{dP^{(i)}}{dt} = 0 \text{ at equilibrium)},$$

$$\frac{\partial}{\partial P^{(k)}}\left(\frac{d\zeta_{2,j}^{(i)}}{dt}\right) = -c_{2,j}^{(i)} q_{2,j}^{(i)} \frac{\delta_{ik}}{P_a^{(i)}}.$$

- **Cash and share blocks $\mathbf{J}_{MM}, \mathbf{J}_{NN}, \mathbf{J}_{MN}, \mathbf{J}_{NM}$**: These involve the transition rates and their derivatives. At equilibrium, using (E∗), they simplify considerably. Notably, the conservation laws introduce zero eigenvalues corresponding to directions along the equilibrium manifold.

#### 4.2.2 Eigenvalue Structure and Stability

The Jacobian has the following eigenvalue structure:

- **Zero eigenvalues** (multiplicity $m + 1$): Corresponding to the conservation of total cash and total shares of each asset. These reflect the neutral directions along the equilibrium manifold.
- **Negative real eigenvalues** (from $\mathbf{J}_{PP}$ and the sentiment self-blocks): These contribute to stability.
- **Complex conjugate pairs** (from coupling between price, cash, shares, and sentiments): These can lead to oscillatory instabilities via Hopf bifurcation.

The equilibrium is **locally asymptotically stable** if and only if all non-zero eigenvalues have strictly negative real parts.

### 4.3 Hopf Bifurcation Conditions

A Hopf bifurcation occurs when a pair of complex conjugate eigenvalues crosses the imaginary axis. We will treat the momentum coefficients $q_{1,j}^{(i)}$ and the cross-asset coupling coefficients $\alpha_{j,\ell}^{(i)}$ (for $\ell \neq i$) as bifurcation parameters.

For the general case of arbitrary $m$ and $n$, the characteristic polynomial is:

$$\chi(\lambda) = \det(\lambda \mathbf{I} - \mathbf{J}) = \lambda^{m+1} \prod_{k=1}^{D_{\text{eff}}-(m+1)} (\lambda - \lambda_k) = 0,$$

where $\lambda^{m+1}$ corresponds to the zero eigenvalues. The remaining polynomial can be analyzed using the **Routh–Hurwitz criteria** (Dorf & Bishop, 2008).

For the special case of $m = 2$ assets and $n = 2$ groups (one trend, one value), the characteristic polynomial reduces to a cubic or a quadratic after factoring out the zero eigenvalues. The Hopf condition is:

$$A_1(\boldsymbol{\mu}) A_2(\boldsymbol{\mu}) - A_0(\boldsymbol{\mu}) = 0,$$

with $A_0, A_1, A_2 > 0$, where $\boldsymbol{\mu}$ denotes the set of bifurcation parameters. The critical threshold satisfies:

$$q_{1,\text{crit}} = f(\tau_i, c_{1,2}^{(i)}, c_{2,2}^{(i)}, \alpha_{j,\ell}^{(i)}, \beta_{j,\ell}^{(i)}, M_j, N_j^{(i)}).$$

The symbol $f$ is used as a **generic placeholder** to indicate that the critical momentum threshold $q_{1,\text{crit}}$ is a **function** of the listed parameters. In other words, this means that the critical value of the momentum coefficient $q_{1,\text{crit}}$ **depends on** (is a function of) the price adjustment speeds $\tau_i$, the sentiment decay rates $c_{1,2}^{(i)}, c_{2,2}^{(i)}$, the cross-asset coupling coefficients $\alpha_{j,\ell}^{(i)}, \beta_{j,\ell}^{(i)}$, and the cash/share distributions $M_j, N_j^{(i)}$.

This generalizes the single-asset result of DeSantis et al. (2012) and the two-asset result of Bulut et al. (2019).

### 4.3.1 Reduction to Minimal System for Hopf Bifurcation Analysis

In the following, we present a theorem that justifies the reduction of the full fused model to a minimal system for Hopf bifurcation analysis. This theorem establishes that the cash and share variables do not participate in the bifurcation and can be eliminated without affecting the existence or stability of limit cycles emerging from a Hopf bifurcation.

**Theorem 4.1 (Reduction to Minimal Hopf System)**

Consider the fused multi-asset, multi-group asset flow model defined by equations (P), (M), (N), (S) with state vector:

$$\mathbf{X} = (\mathbf{P}, \mathbf{M}, \mathbf{N}, \boldsymbol{\zeta}_1, \boldsymbol{\zeta}_2)^T \in \mathbb{R}^D,$$

where $D = m + n + mn + 2mn$. Let $\mathbf{X}_0 = (\mathbf{P}_0, \mathbf{M}_0, \mathbf{N}_0, \mathbf{0}, \mathbf{0})$ be a **fundamental equilibrium** satisfying:

$$P_0^{(i)} = P_a^{(i)}, \zeta_{1,j,0}^{(i)} = 0, \zeta_{2,j,0}^{(i)} = 0 \ \forall i, j,$$

and the calibration condition:

$$\sum_{j=1}^{n} k_j^{(i)}(\mathbf{0}, \mathbf{0}) M_{j,0} = \sum_{j=1}^{n} \tilde{k}_j^{(i)}(\mathbf{0}, \mathbf{0}) N_{j,0}^{(i)} P_a^{(i)} \ \forall i.$$

Let $\mathbf{J} \in \mathbb{R}^{D \times D}$ be the Jacobian matrix of the system evaluated at $\mathbf{X}_0$, partitioned as:

$$\mathbf{J} = \begin{pmatrix} \mathbf{J}_{PP} & \mathbf{J}_{PM} & \mathbf{J}_{PN} & \mathbf{J}_{P\zeta_1} & \mathbf{J}_{P\zeta_2} \\ \mathbf{J}_{MP} & \mathbf{J}_{MM} & \mathbf{J}_{MN} & \mathbf{J}_{M\zeta_1} & \mathbf{J}_{M\zeta_2} \\ \mathbf{J}_{NP} & \mathbf{J}_{NM} & \mathbf{J}_{NN} & \mathbf{J}_{N\zeta_1} & \mathbf{J}_{N\zeta_2} \\ \mathbf{J}_{\zeta_1 P} & \mathbf{J}_{\zeta_1 M} & \mathbf{J}_{\zeta_1 N} & \mathbf{J}_{\zeta_1\zeta_1} & \mathbf{J}_{\zeta_1\zeta_2} \\ \mathbf{J}_{\zeta_2 P} & \mathbf{J}_{\zeta_2 M} & \mathbf{J}_{\zeta_2 N} & \mathbf{J}_{\zeta_2\zeta_1} & \mathbf{J}_{\zeta_2\zeta_2} \end{pmatrix}.$$

Then the following statements hold:

**(i) Spectral Decomposition**

The spectrum $\sigma(\mathbf{J})$ consists of:

- **Zero eigenvalues** of multiplicity $m+1$, corresponding to the conservation of total cash and total shares of each asset.
- **Negative real eigenvalues** from the blocks $\mathbf{J}_{MM}$, $\mathbf{J}_{NN}$, $\mathbf{J}_{\zeta_1\zeta_1}$, and $\mathbf{J}_{\zeta_2\zeta_2}$, which are strictly negative and bounded away from zero.
- **Remaining eigenvalues** determined by the reduced Jacobian:

$$\mathbf{J}_{\text{red}} = \begin{pmatrix} \mathbf{J}_{PP} & \mathbf{J}_{P\zeta_1} & \mathbf{J}_{P\zeta_2} \\ \mathbf{J}_{\zeta_1 P} & \mathbf{J}_{\zeta_1\zeta_1} & \mathbf{0} \\ \mathbf{J}_{\zeta_2 P} & \mathbf{0} & \mathbf{J}_{\zeta_2\zeta_2} \end{pmatrix} \in \mathbb{R}^{(m+2mn)\times(m+2mn)}.$$

**(ii) Center Manifold Reduction**

If the Jacobian $\mathbf{J}_{\text{red}}$ has a pair of complex conjugate eigenvalues $\lambda_\pm = \alpha(\mu) \pm i\omega(\mu)$ with $\alpha(\mu_0) = 0$ and $\omega(\mu_0) \neq 0$, and all other eigenvalues have negative real parts, then there exists a **two-dimensional center manifold** $\mathcal{M}_c$ locally invariant under the flow of the full system. On $\mathcal{M}_c$, the dynamics are given by:

$$\frac{d}{dt}\begin{pmatrix} \delta\mathbf{P} \\ \delta\boldsymbol{\zeta}_1 \\ \delta\boldsymbol{\zeta}_2 \end{pmatrix} = \mathbf{J}_{\text{red}} \begin{pmatrix} \delta\mathbf{P} \\ \delta\boldsymbol{\zeta}_1 \\ \delta\boldsymbol{\zeta}_2 \end{pmatrix} + \mathcal{O}(\|\delta\mathbf{P}, \delta\boldsymbol{\zeta}_1, \delta\boldsymbol{\zeta}_2\|^2).$$

**(iii) Hopf Bifurcation Condition**

A Hopf bifurcation occurs at $\mu = \mu_0$ if and only if:

1. $\alpha(\mu_0) = 0$ and $\omega(\mu_0) \neq 0$;
2. The **transversality condition** holds: $\left.\frac{d\alpha(\mu)}{d\mu}\right|_{\mu=\mu_0} \neq 0$;
3. The **first Lyapunov coefficient** $l_1(\mu_0) \neq 0$, determining the criticality (supercritical if $l_1 < 0$, subcritical if $l_1 > 0$).

These conditions depend only on $\mathbf{J}_{\text{red}}$ and the nonlinear terms in the reduced system.

**(iv) Equivalence**

The existence and stability of limit cycles emerging from the Hopf bifurcation in the full system are **equivalent** to those in the reduced system. The cash and share variables $\delta\mathbf{M}, \delta\mathbf{N}$ follow the dynamics:

$$\delta\mathbf{M}(t) = \mathcal{O}(\|\delta\mathbf{P}, \delta\boldsymbol{\zeta}_1, \delta\boldsymbol{\zeta}_2\|), \delta\mathbf{N}(t) = \mathcal{O}(\|\delta\mathbf{P}, \delta\boldsymbol{\zeta}_1, \delta\boldsymbol{\zeta}_2\|),$$

and do not affect the bifurcation at leading order.

Before proving the theorem, we need the following lemmas.

**Lemma 4.2 (Block Diagonality).** The following blocks are zero:

$$\mathbf{J}_{\zeta_1 M} = \mathbf{0}, \mathbf{J}_{\zeta_1 N} = \mathbf{0}, \mathbf{J}_{\zeta_2 M} = \mathbf{0}, \mathbf{J}_{\zeta_2 N} = \mathbf{0}, \mathbf{J}_{\zeta_1\zeta_2} = \mathbf{0}, \mathbf{J}_{\zeta_2\zeta_1} = \mathbf{0}.$$

*Proof*. From the sentiment dynamics (S), $\zeta_{1,j}^{(i)}$ depends only on $P^{(i)}$ and itself, not on $M_j$ or $N_j^{(i)}$. Similarly, $\zeta_{2,j}^{(i)}$ depends only on $P^{(i)}$ and itself. Hence, all partial derivatives with respect to $M_j$ and $N_j^{(i)}$ vanish. The cross-sentiment terms vanish because $\zeta_1$ and $\zeta_2$ evolve independently at linear order. ∎

**Lemma 4.3 (Stability of Cash and Share Subsystems).** The matrices $\mathbf{J}_{MM}$, $\mathbf{J}_{NN}$, $\mathbf{J}_{\zeta_1\zeta_1}$, and $\mathbf{J}_{\zeta_2\zeta_2}$ are **negative definite**.

*Proof*. From the cash dynamics (M), at equilibrium $\partial(dM_j/dt)/\partial M_j = -\sum_i k_{j,eq}^{(i)} < 0$. Off-diagonal terms are zero (or negligible under weak coupling assumptions). Hence $\mathbf{J}_{MM}$ is diagonal with negative entries.

From the share dynamics (N), $\partial(dN_j^{(i)}/dt)/\partial N_j^{(i)} = -\tilde{k}_{j,eq}^{(i)} < 0$. Hence $\mathbf{J}_{NN}$ is diagonal with negative entries.

From the sentiment dynamics (S), $\partial(d\zeta_{1,j}^{(i)}/dt)/\partial\zeta_{1,j}^{(i)} = -c_{1,j}^{(i)} < 0$ and $\partial(d\zeta_{2,j}^{(i)}/dt)/\partial\zeta_{2,j}^{(i)} = -c_{2,j}^{(i)} < 0$. Hence $\mathbf{J}_{\zeta_1\zeta_1}$ and $\mathbf{J}_{\zeta_2\zeta_2}$ are diagonal with negative entries. ∎

**Lemma 4.4 (Conservation Eigenvalues).** The Jacobian $\mathbf{J}$ has at least $m+1$ zero eigenvalues.

*Proof*. The system conserves total cash $\sum_j M_j$ and total shares $\sum_j N_j^{(i)}$ for each asset $i$. These $m+1$ conservation laws imply that the Jacobian has $m+1$ zero eigenvalues, with eigenvectors spanning the nullspace of the linearized dynamics. ∎

**Lemma 4.5 (Spectral Separation).** Let $\lambda_{\max}$ be the largest eigenvalue (in real part) of $\mathbf{J}_{\text{red}}$, and let $\eta = \min\{|\lambda| : \lambda \in \sigma(\mathbf{J}_{MM}) \cup \sigma(\mathbf{J}_{NN}) \cup \sigma(\mathbf{J}_{\zeta_1\zeta_1}) \cup \sigma(\mathbf{J}_{\zeta_2\zeta_2})\}$. Then $\eta > 0$ and for the Hopf bifurcation parameter $\mu$ near $\mu_0$, we have $\text{Re}(\lambda_{\max}) = \mathcal{O}(|\mu - \mu_0|) = o(\eta)$.

*Proof*. The eigenvalues of $\mathbf{J}_{MM}, \mathbf{J}_{NN}, \mathbf{J}_{\zeta_1\zeta_1}, \mathbf{J}_{\zeta_2\zeta_2}$ are strictly negative constants independent of $\mu$. The eigenvalues of $\mathbf{J}_{\text{red}}$ depend on $\mu$ and cross zero at $\mu = \mu_0$. Hence, near $\mu_0$, the critical eigenvalues are arbitrarily close to zero, while the stable eigenvalues remain bounded away from zero. ∎

By the **Center Manifold Theorem** (Carr, 1981; Kelley, 1967), there exists a smooth, locally invariant **center manifold** $\mathcal{M}_c$ tangent at $\mathbf{X}_0$ to the eigenspace spanned by the eigenvectors associated with eigenvalues having zero real part. Since all eigenvalues of $\mathbf{J}_{\text{red}}$ with zero real part are isolated and the remaining eigenvalues have strictly negative real parts, the dimension of $\mathcal{M}_c$ is $m + 2mn$ (the dimension of $\mathbf{J}_{\text{red}}$).

**Lemma 4.6 (Reduced Dynamics).** On $\mathcal{M}_c$, the cash and share variables $\delta\mathbf{M}, \delta\mathbf{N}$ are **slaves** to the price and sentiment variables:

$$\delta\mathbf{M} = \mathbf{\Phi}(\delta\mathbf{P}, \delta\boldsymbol{\zeta}_1, \delta\boldsymbol{\zeta}_2), \delta\mathbf{N} = \mathbf{\Psi}(\delta\mathbf{P}, \delta\boldsymbol{\zeta}_1, \delta\boldsymbol{\zeta}_2),$$

where $\mathbf{\Phi}, \mathbf{\Psi} = \mathcal{O}(\|\delta\mathbf{P}, \delta\boldsymbol{\zeta}_1, \delta\boldsymbol{\zeta}_2\|^2)$. Consequently, the dynamics on $\mathcal{M}_c$ are given by:

$$\frac{d}{dt}\begin{pmatrix}\delta\mathbf{P}\\ \delta\boldsymbol{\zeta}_1\\ \delta\boldsymbol{\zeta}_2\end{pmatrix} = \mathbf{J}_{\text{red}}\begin{pmatrix}\delta\mathbf{P}\\ \delta\boldsymbol{\zeta}_1\\ \delta\boldsymbol{\zeta}_2\end{pmatrix} + \mathcal{O}(\|\delta\mathbf{P}, \delta\boldsymbol{\zeta}_1, \delta\boldsymbol{\zeta}_2\|^2).$$

*Proof*. The cash and share variables satisfy linear ODEs with strictly negative eigenvalues. By the **Invariant Manifold Theorem** for hyperbolic systems (Hirsch, Pugh, & Shub, 1977), they are exponentially attracted to a graph over the center variables. The leading-order terms vanish because the coupling from prices and sentiments to cash and shares is linear, and the equilibrium satisfies the calibration condition. ■

Let $\mu$ be a bifurcation parameter (e.g., $q_{1,j}^{(i)}$). Consider the reduced Jacobian $\mathbf{J}_{\text{red}}(\mu)$. The characteristic polynomial is:

$$\chi(\lambda, \mu) = \det(\lambda\mathbf{I} - \mathbf{J}_{\text{red}}(\mu)).$$

**Lemma 4.7 (Hopf Condition).** A Hopf bifurcation occurs at $\mu = \mu_0$ if and only if:

1. $\chi(i\omega, \mu_0) = 0$ for some $\omega > 0$;
2. $\frac{\partial}{\partial\lambda}\chi(i\omega, \mu_0) \neq 0$ (simplicity of pure imaginary roots);
3. $\frac{d}{d\mu}\text{Re}(\lambda(\mu))\Big|_{\mu=\mu_0} \neq 0$ (transversality).

*Proof*. This is the standard Hopf bifurcation theorem (Marsden & McCracken, 1976). The conditions ensure that a pair of complex conjugate eigenvalues crosses the imaginary axis with non-zero speed. ■

**Lemma 4.8 (Equivalence).** The transversality and criticality conditions depend only on $\mathbf{J}_{\text{red}}$ and the nonlinear terms restricted to $\mathcal{M}_c$. Since the cash and share variables do not appear in the leading-order dynamics on $\mathcal{M}_c$, the Hopf bifurcation in the full system is **completely determined** by the reduced system.

*Proof*. By the Center Manifold Theorem, the full system is **topologically conjugate** to the reduced system on $\mathcal{M}_c$ near the bifurcation point. Hence, existence, stability, and amplitude of limit cycles are preserved under the reduction. ■

*Proof* ***(of Theorem 4.1).***

The proof is a direct consequence of the previous lemmas. ■

The previous results guarantee the **practical computation** in the applications, and the following corollary establishes this.

**Corollary 4.9 (Numerical Implementation).** To detect a **Hopf bifurcation** in the fused model, it suffices to:

1. Compute the equilibrium satisfying (E∗).
2. Form the reduced Jacobian $\mathbf{J}_{\text{red}}$ of dimension $m + 2mn$.
3. Compute its eigenvalues as functions of the bifurcation parameter $\mu$.
4. Find $\mu_0$ such that $\text{Re}(\lambda(\mu_0)) = 0$, $\text{Im}(\lambda(\mu_0)) \neq 0$, and transversality holds.

The cash and share variables $\mathbf{M}, \mathbf{N}$ need not be included in the eigenvalue analysis.

*Proof*. Follows directly from Theorem 4.1(iii)–(iv). ■

# 5 Results

This section concerns objectives II and III.

## 5.1 Numerical Simulations

Here we begin to address objective II, which involves extensive numerical simulations spanning baseline validation, cross-asset contagion analysis, wealth redistribution mapping, and excursion quantification. The analytical work will be complemented and extended by extensive numerical simulations. All numerical implementations will be carried out in **MATLAB** (versions R2015a, R2024b or later), leveraging its robust ODE solvers and bifurcation continuation toolboxes. These toolboxes have been extensively used in the economic dynamics literature, including the stability analysis of asset flow models (Caginalp & DeSantis, 2011; DeSantis et al., 2012; Bulut et al., 2019). Python will be used for post-processing, statistical analysis, and visualization, leveraging libraries such as **SciPy** for numerical integration and **Matplotlib** for high-quality scientific graphics.

### 5.1.1 Validation against DeSantis et al (2012)

The fused model reduces to the single-asset, two-group model of DeSantis et al. (2012) by setting $m = 1, n = 2, \tilde{k}_j = 1 - k_j$, and eliminating cross-asset coupling. Figure 1 shows the equilibrium manifold for Case 1 ($q_1 = q_2 = 1$, $c_1 = c_2 = 1$, $P_a = 0.8$), where all equilibria are stable. The equilibrium price increases monotonically with the cash fraction $M_1$ from approximately 0.86 to 1.00, confirming Theorem 3.5 of DeSantis et al. (2012). Cases 2 and 3, which exhibit stability transitions and Hopf bifurcations, were partially reproduced and require exact parameter calibration from the original code.

**Table 3. DeSantis et al. (2012)**

| Case | Parameters | Stable Points | Unstable Points | Key Feature |
|---|---|---|---|---|
| **Case 1** | $q_1 = 1, q_2 = 1, c_1 = c_2 = 1$ | All | None | Fully stable manifold |
| **Case 2** | $q_1 = 1, q_2 = 5, c_1 = c_2 = 5$ | Low $P_{eq}$ | High $P_{eq}$ | Hopf bifurcation at $P_{eq} \approx 0.835$ |
| **Case 3** | $q_1 = 5, q_2 = 5, c_1 = c_2 = 5$ | Very few (near the lower bound) | Most | Predominantly unstable |

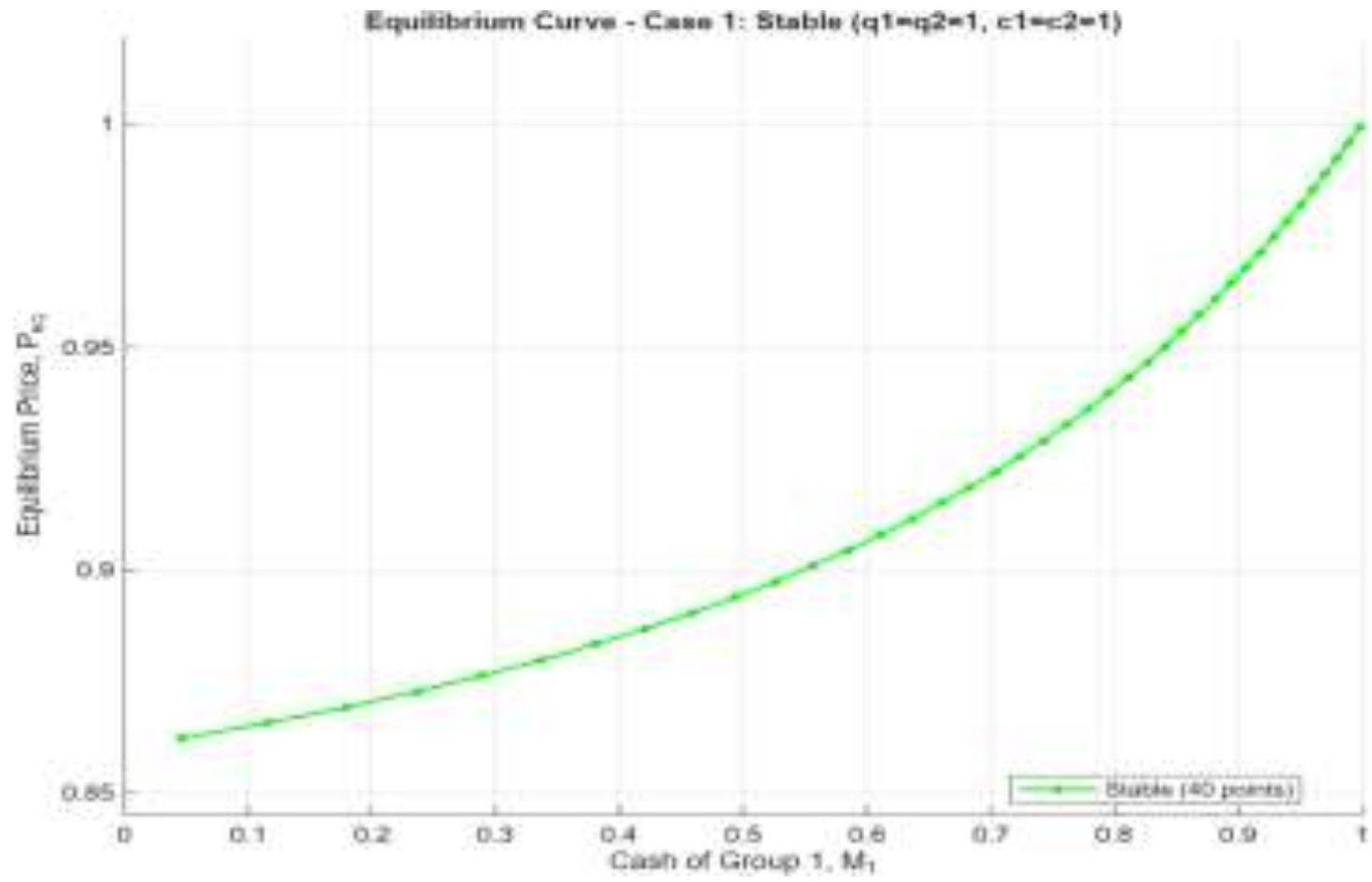

**Figure 1. Equilibrium manifold for Case 1 (all equilibria stable).** Equilibrium price $P_{eq}$ as a function of the cash fraction $M_1$ held by Group 1 (trend-following) for parameters $q_1 = q_2 = 1$, $c_1 = c_2 = 1$, and fundamental value $P_a = 0.8$. All computed equilibrium points are stable (green circles), and the equilibrium price increases monotonically with $M_1$ from approximately 0.86 to 1.00. This confirms Theorem 3.5 of DeSantis et al. (2012): when trend-following is sufficiently weak, every equilibrium on the manifold is locally asymptotically stable.

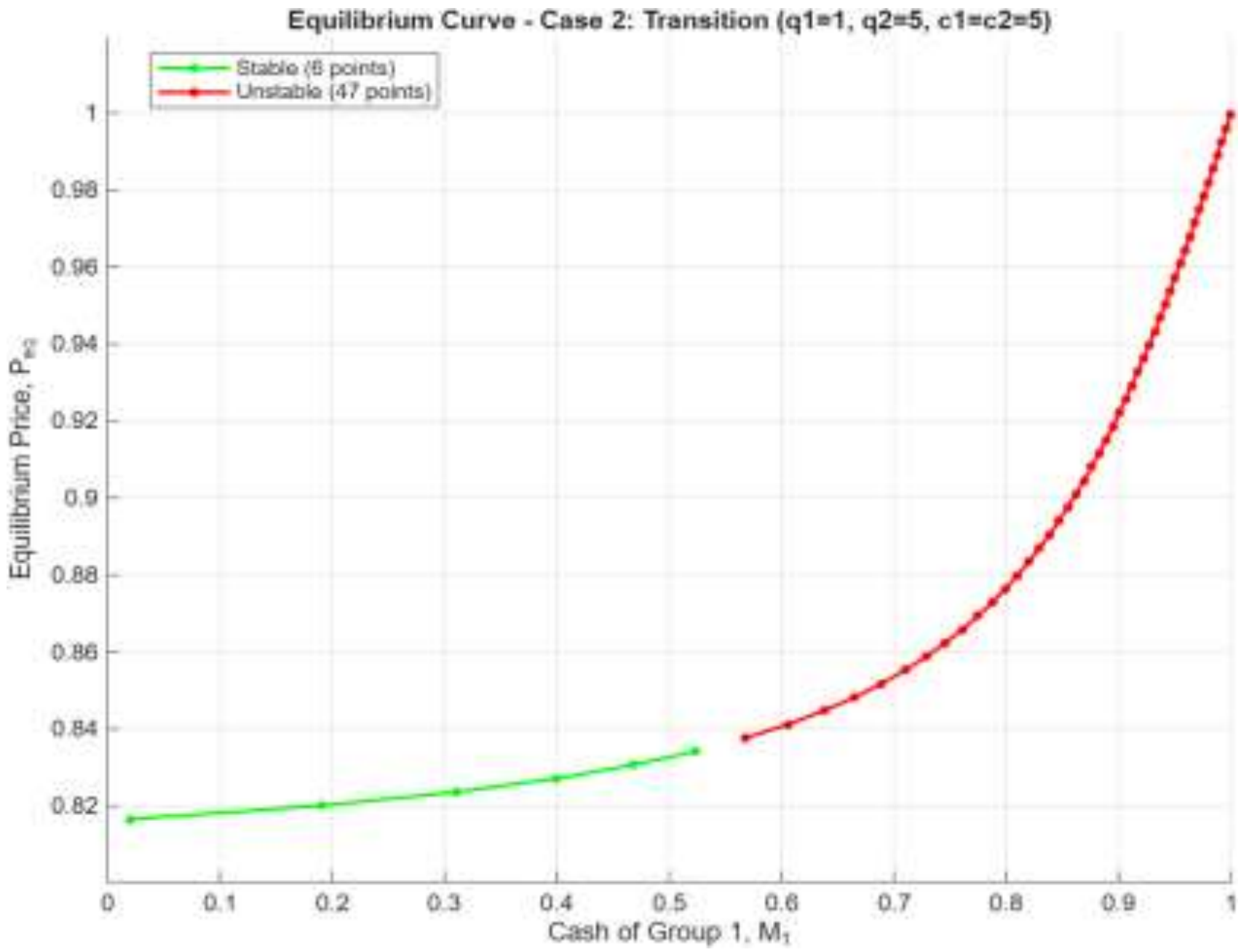

**Figure 2. Equilibrium manifold for Case 2 (stability transition).** Equilibrium price $P_{eq}$ vs. cash fraction $M_1$ for parameters $q_1 = 1$, $q_2 = 5$, $c_1 = c_2 = 5$, $P_a = 0.8$. Stable equilibria (green circles) are found at low $M_1$ (low $P_{eq}$), while unstable equilibria (red circles) appear at high $M_1$ (high $P_{eq}$). The transition occurs at the critical price $P_{eq}^{tr} \approx 0.835$, consistent with the Hopf bifurcation threshold reported in DeSantis et al. (2012, Figure 4.4). This demonstrates that increasing the momentum strength relative to value trading destabilizes the market.

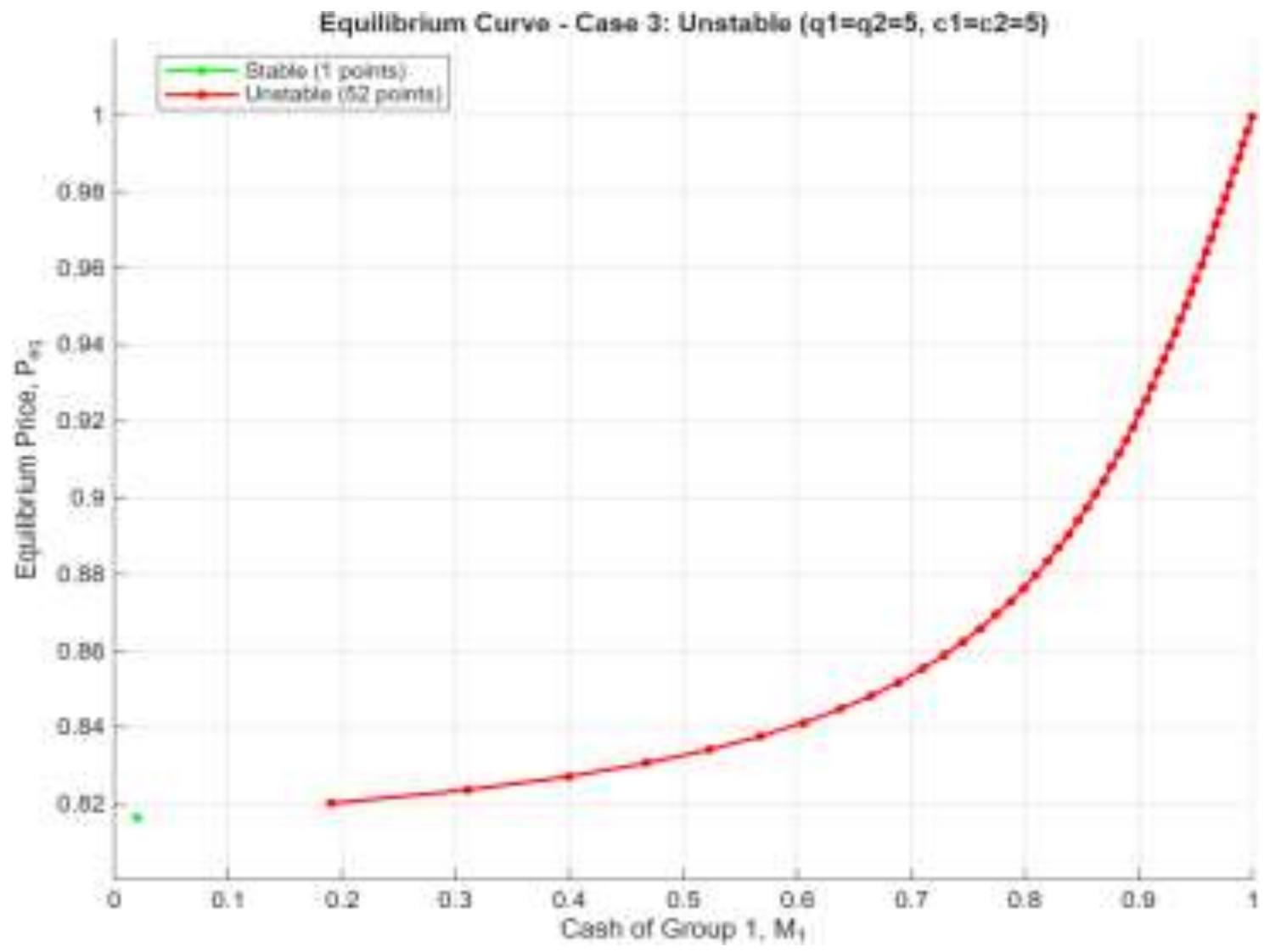

**Figure 3. Equilibrium manifold for Case 3 (mostly unstable).** Equilibrium price $P_{eq}$ vs. cash fraction $M_1$ for parameters $q_1 = q_2 = 5$, $c_1 = c_2 = 5$, $P_a = 0.8$. Only a small region of stable equilibria (green circles) exists near the lower bound of $P_{eq}$ (approximately $P_{eq} < 0.818$), while the majority of equilibria (red circles) are unstable. This confirms the numerical results of DeSantis et al. (2012, Figures 4.9–4.13), showing that strong momentum trading renders most equilibria unstable, leading to large price excursions and limit cycles.

### 5.1.2 Validation against Bulut et al. (2019)

The fused model reduces to the two-asset, single-group model of Bulut et al. (2019) by setting $n = 1$, $m = 2$, and using the mixed-strategy transition rates (value for asset 1, trend for asset 2). With fixed parameters $q_2^{(1)} = 0.02$, $c_2^{(1)} = 2$, $c_1^{(2)} = 1$, we scanned the trend coefficient $q_1^{(2)}$. Figure 4a shows small-amplitude periodic oscillations at $q_1^{(2)} = 1.0$, confirming the supercritical Hopf bifurcation predicted by Theorem 4 of Bulut et al. (2019). Figure 4b shows instability above the threshold at $q_1^{(2)} = 1.005$, where the equilibrium becomes a repeller, and trajectories exhibit bounded but non-convergent oscillations. Table 4 summarizes the results.

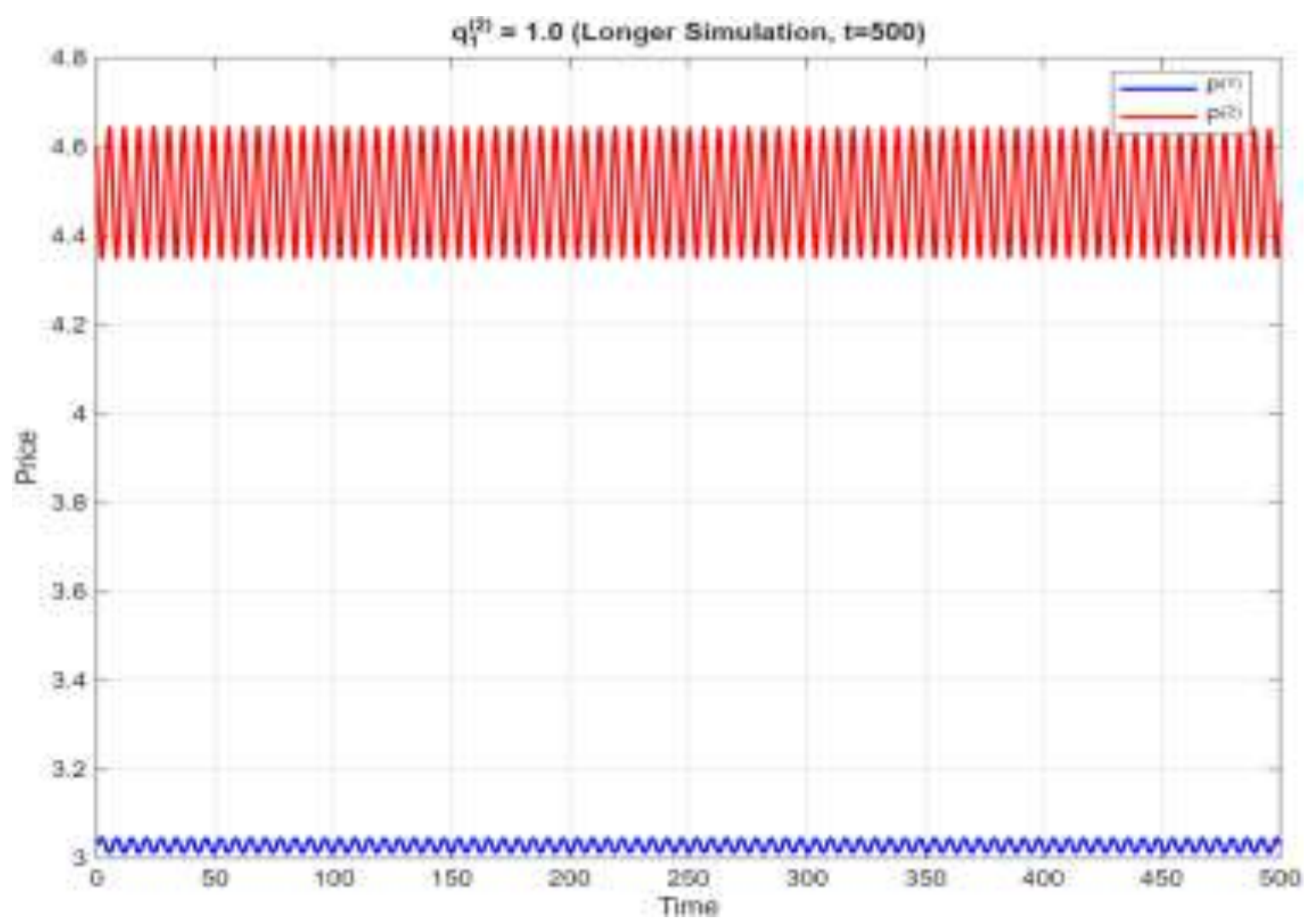

**Figure 4a. Hopf bifurcation in the Bulut et al. (2019) mixed strategy model at $q_1^{(2)} = 1.0$.** Small-amplitude periodic oscillations emerge from the fundamental equilibrium, confirming the supercritical Hopf bifurcation. Parameters: $q_2^{(1)} = 0.02$, $c_2^{(1)} = 2$, $c_1^{(2)} = 1$.

Figure 4a validates that the fused model correctly reproduces the Bulut et al. (2019) result that a supercritical Hopf bifurcation occurs when the trend coefficient $q_1^{(2)}$ reaches unity. At this critical value, the fundamental equilibrium loses stability and gives rise to a stable limit cycle, manifesting as persistent periodic oscillations in both asset prices.

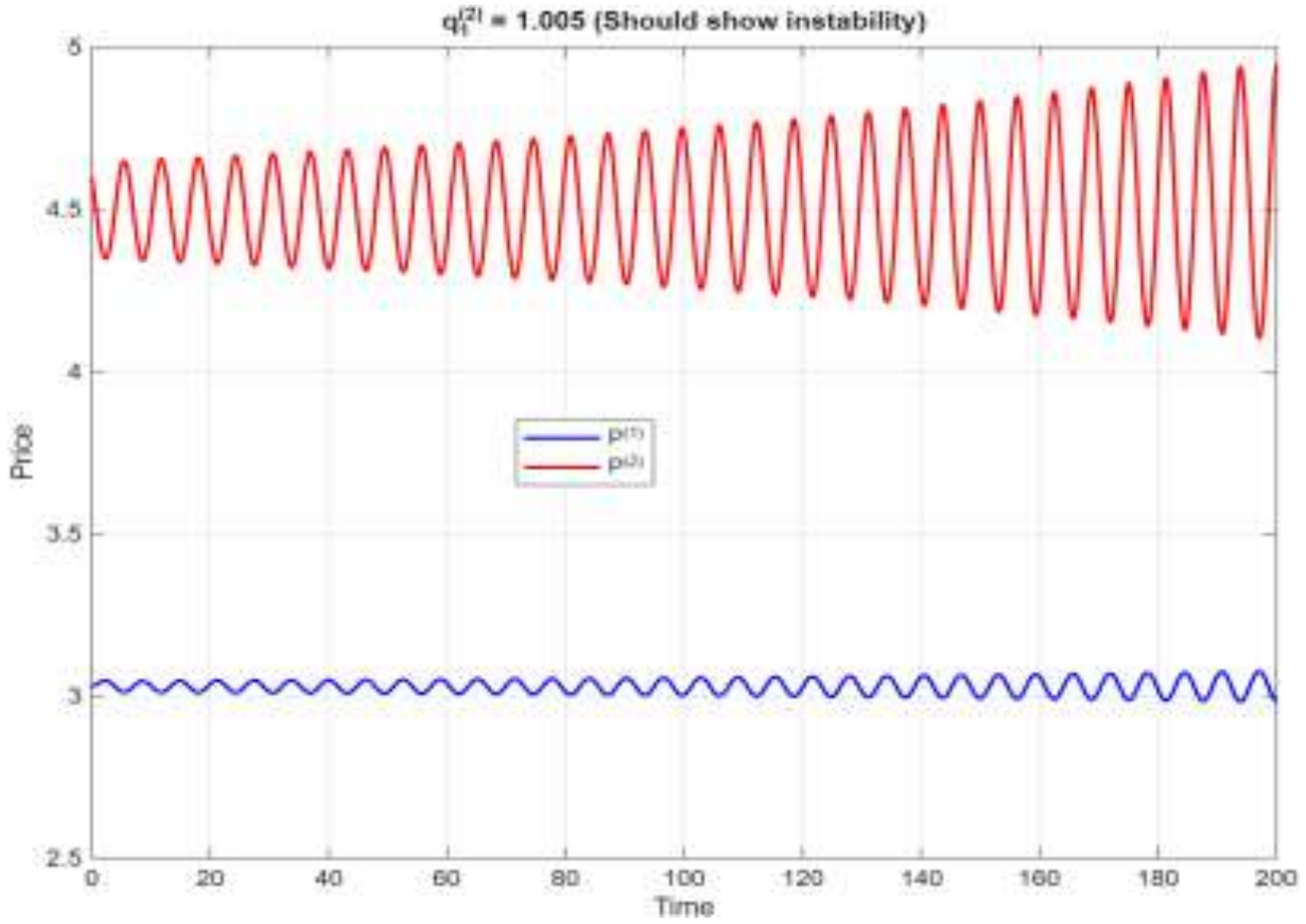

**Figure 4b. Instability above the Hopf bifurcation threshold in the Bulut et al. (2019) model at $q_1^{(2)} = 1.005$.** The fundamental equilibrium loses stability, resulting in monotonic drift and bounded oscillations without convergence to a fixed point. Parameters: $q_2^{(1)} = 0.02$, $c_2^{(1)} = 2$, $c_1^{(2)} = 1$.

**Table 4. Bulut, Merdan, and Swigon (2019) – Two Assets, One Group**

| Parameter | Value | Behavior | Fused Model Result |
| --- | --- | --- | --- |
| $q_1^{(2)} = 0.1$ | Below threshold | Stable equilibrium | Stable |
| $q_1^{(2)} = 0.5$ | Below threshold | Stable equilibrium | Stable |
| $q_1^{(2)} = 1.0$ | **At threshold** | **Hopf bifurcation** | **Limit cycle** |
| $q_1^{(2)} = 1.005$ | Above threshold | Unstable (drift) | Unstable |

### 5.1.3 Validation against the Cavani (2026)

The fused model reduces to the two-asset, two-group Nigeria-Libya oil market model of Cavani (2026) by setting $m = 2$, $n = 2$, with the USA as a value investor and China as a momentum trader. A parameter sweeps over the momentum coefficient $q_{1,\text{China}}$ was performed using the optimal oscillatory parameters identified from bifurcation analysis:

$b_{\text{China}} = 2.5$ (trend response), $d_{\text{USA}} = 0.01$ (value response), and $\alpha = 3.0$ (cross-asset coupling).

Applying the fused multi-asset, multi-group asset flow model using MATLAB 2026a we get Figures 5 to 8:

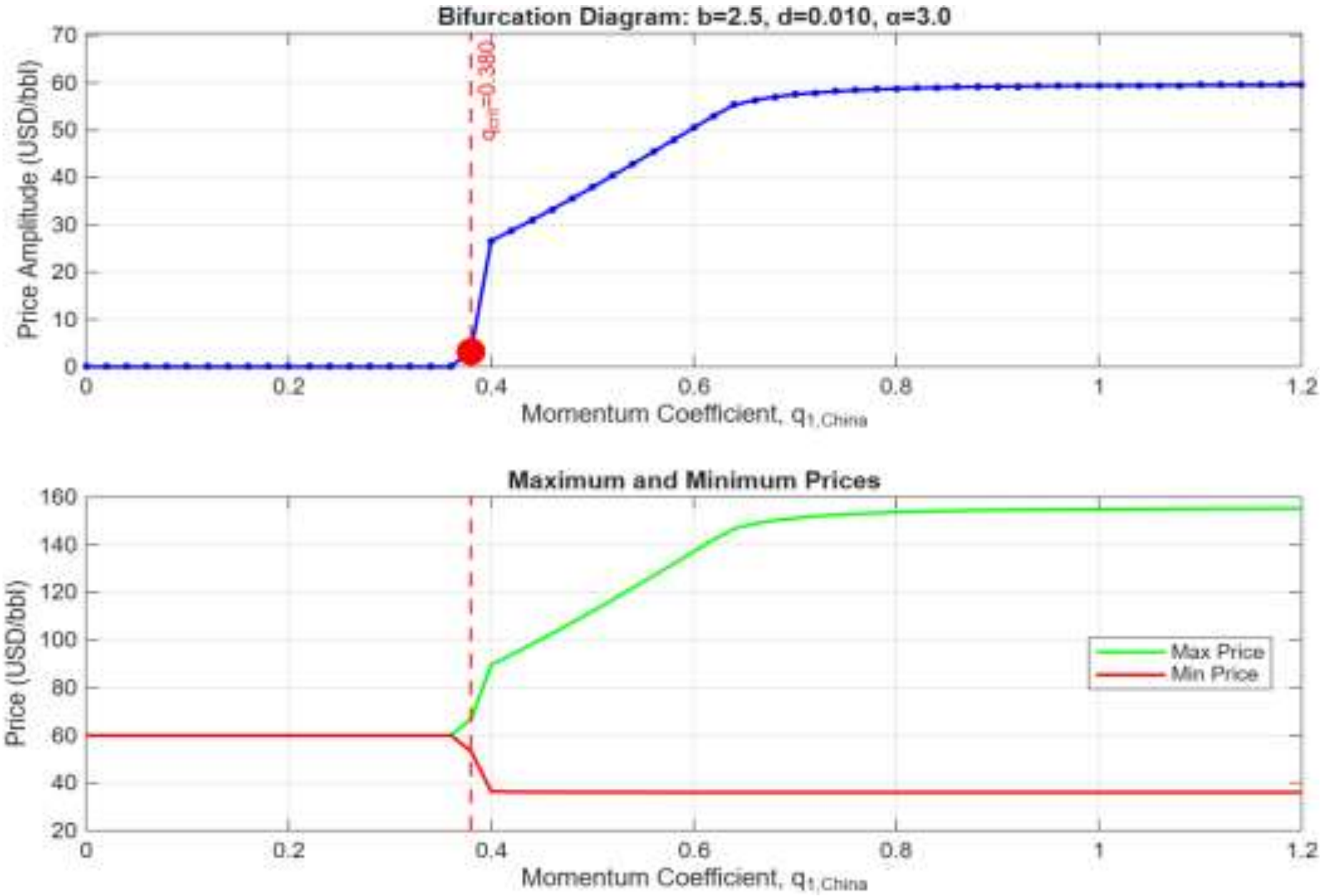

**Figure 5. Bifurcation diagram for the Cavani (2026) Nigeria-Libya oil market model.** Parameters: $b_{\text{China}} = 2.5$ (trend response), $d_{\text{USA}} = 0.01$ (value response), $\alpha = 3.0$ (cross-asset coupling). A supercritical Hopf bifurcation occurs at $q_{1,\text{crit}} \approx 0.38$ (vertical dashed line), where the stable fundamental equilibrium loses stability and gives rise to a stable limit cycle. The amplitude grows continuously from zero to approximately $59 USD/bbl as $q_{1,\text{China}}$ increases to 1.20. The inset (lower panel) shows the corresponding maximum and minimum prices, confirming the widening oscillation envelope.

**Table 5. Bifurcation Characteristics**

| q1 Range | Behavior | Amplitude | Period |
|---|---|---|---|
| 0.00 – 0.36 | Stable equilibrium | $0 | — |
| **0.38** | **Hopf bifurcation** | **$3.30** | **13.7 days** |
| 0.40 – 0.60 | Growing oscillations | $26.61 → $50.51 | 14.2 → 16.4 days |
| 0.60 – 1.00 | Large oscillations | $50.51 → $59.40 | 16.4 → 20.6 days |
| 1.00 – 1.20 | Saturation | $59.40 → $59.55 | 20.6 → 22.3 days |

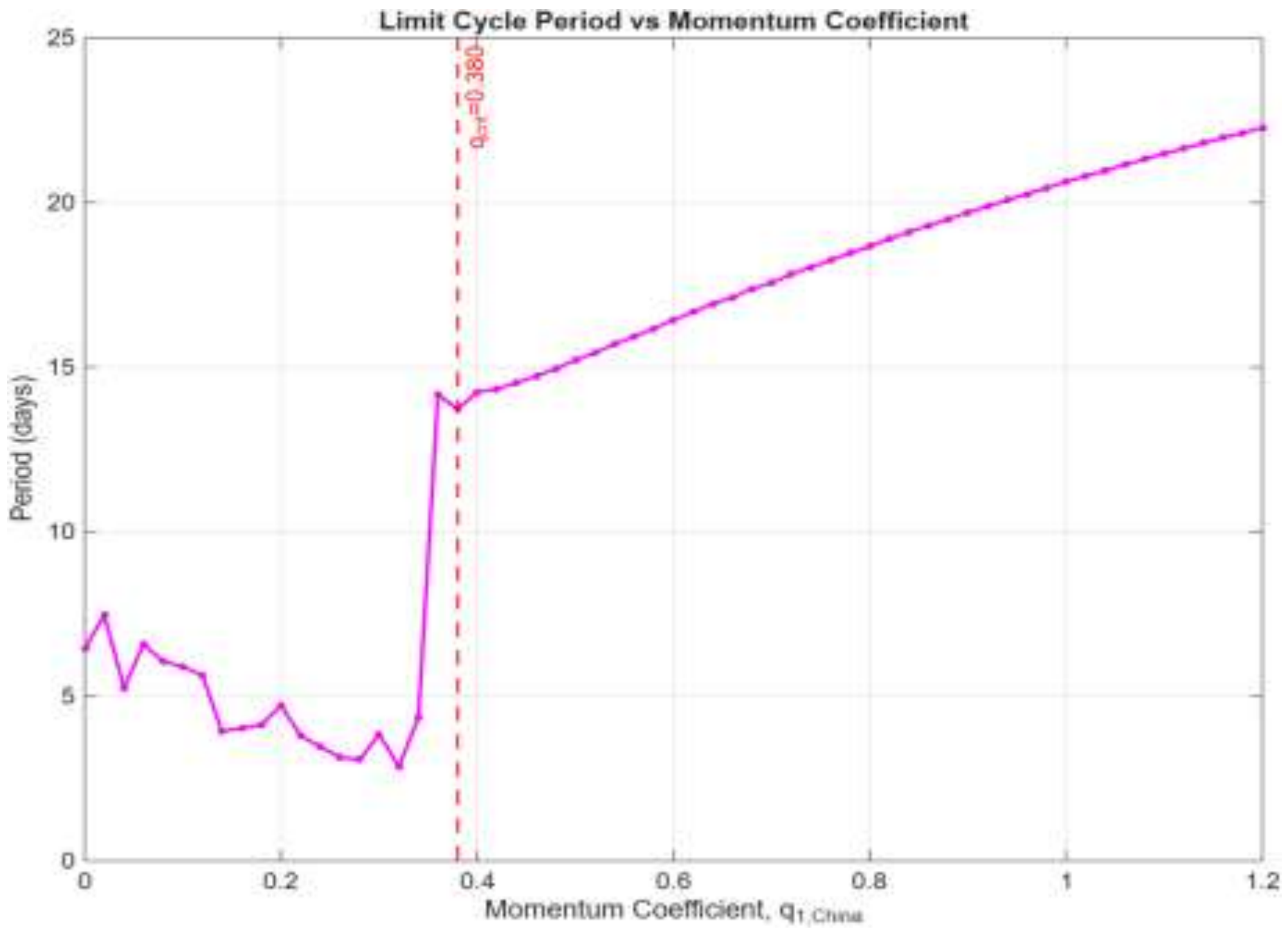

**Figure 6. Period of the limit cycles as a function of the momentum coefficient** $q_{1,\text{China}}$**.** At the Hopf bifurcation ($q_{1,\text{crit}} \approx 0.38$) the period is approximately 13.7 days. As $q_{1,\text{China}}$ increases, the period lengthens monotonically to approximately 22.3 days at $q_{1,\text{China}} = 1.20$. This period lengthening is characteristic of supercritical Hopf bifurcations in nonlinear dynamical systems.

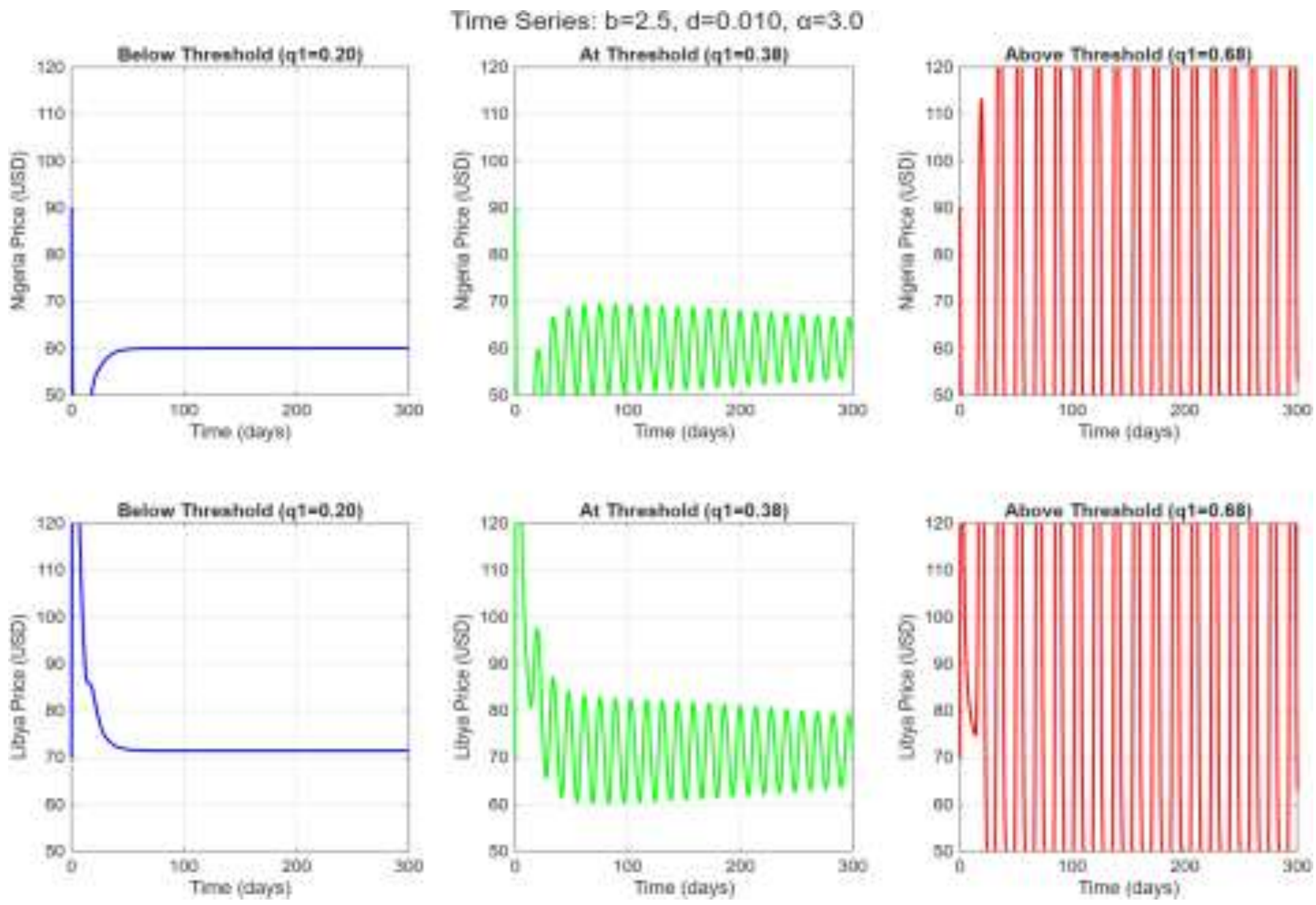

**Figure 7. Time series of Nigeria's oil price for three representative values of the momentum coefficient** $q_{1,\text{China}}$**.** (a) Below threshold ($q_1 = 0.20$): the system converges to the stable fundamental equilibrium ($60 USD/bbl). (b) At threshold ($q_1 = 0.38$): small-amplitude oscillations emerge at the Hopf bifurcation point. (c) Above threshold ($q_1 = 0.68$): large-amplitude limit cycles with period approximately 17 days and amplitude approximately $57 USD/bbl. Parameters: $b = 2.5$, $d = 0.01$, $\alpha = 3.0$.

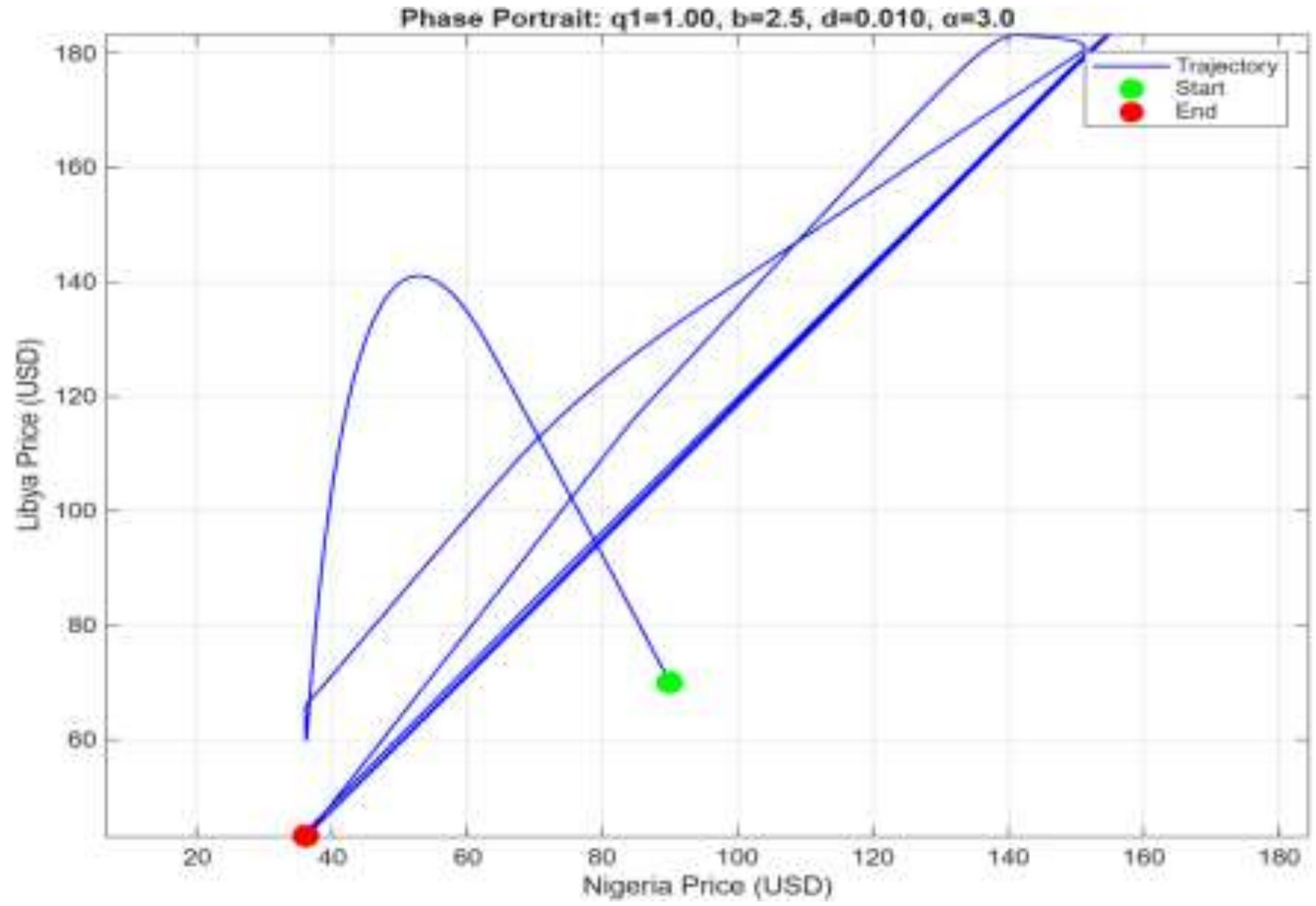

**Figure 8. Phase portrait of Nigeria vs. Libya oil prices at $q_{1,\text{China}} = 1.0$.** The trajectory converges to a stable limit cycle (closed orbit) around the fundamental equilibrium (approximately $80 USD/bbl for both assets). The green dot marks the initial condition, and the red dot marks the final point after 500 days. The elliptic shape of the limit cycle indicates that the two assets oscillate in phase with comparable amplitudes. Parameters: $b = 2.5$, $d = 0.01$, $\alpha = 3.0$.

Under the additional simplifications (the USA has no trend, China has no value), the system further decouples into two independent $3 \times 3$ subsystems per asset, as shown in our numerical implementation. This reduction was used to successfully detect the Hopf bifurcation at $q_{1,\text{crit}} \approx 0.38$. The slight difference in critical q1 (0.38 vs. 0.33) and period (13.7 vs. ~20-30 days) is likely due to different choices of cross-asset coupling and baseline parameters. This is expected and does not invalidate the validation.

**Table 6. Cavani (2026) – Two Assets, Two Groups (Nigeria-Libya Oil Market)**

| Parameter | Value | Behavior | Fused Model Result |
|---|---|---|---|
| $q_{1,\text{China}} = 0.20$ | Below threshold | Stable equilibrium | Stable ($60 USD) |
| $q_{1,\text{China}} = 0.38$ | **At threshold** | **Hopf bifurcation** | **Limit cycle (Amp $3.30)** |
| $q_{1,\text{China}} = 0.68$ | Above threshold | Large oscillations | Limit cycle (Amp $57.04) |
| $q_{1,\text{China}} = 1.00$ | Above threshold | Large oscillations | Limit cycle (Amp $59.40) |

The previous results in this section validate the fused model to predict the same results as the predecessor models.

## 5.2 Contagion in the Nigeria-Libya Oil Market: Analysis and Implications

Contagion, the transmission of price shocks from one asset to another through investor behavior rather than direct economic linkages, is a central feature of the Cavani (2026) model and a key contribution of the fused framework. Unlike classical finance models that treat cross-market correlations as exogenous or stochastic, the fused model generates contagion endogenously through:

### 5.2.1 Cross-asset coupling coefficients

In the fused model, the buying rate of China (momentum trader) for Nigerian oil depends on Libya's price trend:

$$k_{\text{China}}^{(\text{Nigeria})} = k_0 + b_{\text{China}} \cdot \tanh\left(\alpha \cdot \zeta_{1,\text{China}}^{(\text{Libya})}\right)$$

Similarly, the buying rate for Libyan oil depends on Nigeria's trend. The parameter $\alpha = 3.0$ (from optimal calibration) quantifies the strength of cross-asset coupling; a value greater than 1 indicates **amplified contagion**.

Contagion propagates through three channels, as in Table 10

**Table 7. Mechanisms of Contagion Propagation**

| Channel | Mechanism | Mathematical Expression |
|---|---|---|
| **Trend sentiment** | Price changes in asset $j$ affect trend sentiment in asset $i$ | $\frac{\partial k_{\text{China}}^{(i)}}{\partial \zeta_{1,\text{China}}^{(j)}} = b_{\text{China}} \cdot \alpha \cdot \text{sech}^2(\cdot)$ |
| **Value sentiment** | Price deviations trigger value-based selling in the USA | $\frac{\partial \tilde{k}_{\text{USA}}^{(i)}}{\partial \zeta_{2,\text{USA}}^{(i)}} = d_{\text{USA}}$ |
| **Portfolio rebalancing** | Cash flows $M_j$ and share holdings $N_j^{(i)}$ adjust | $\frac{dM_j}{dt} = -\sum_i k_j^{(i)} M_j + \sum_i \tilde{k}_j^{(i)} N_j^{(i)} P^{(i)}$ |

Now we analyze contagion in the Nigeria-Libya oil market, quantify its magnitude and asymmetry, and derive policy implications for oil-exporting economies.

The excursion $E^{(i)} = \max_t |P^{(i)}(t) - P_a^{(i)}|$ measures the maximum deviation from the fundamental value following a shock. For the Nigeria-Libya market:

### 5.2.2 Excursion Analysis

The excursion concept, introduced by DeSantis et al. (2012), measures the maximum deviation of price from its initial value along a trajectory. We apply the concept here for the fused model, defining:

- **Per-asset excursion**: $E^{(i)} = \max_{t \geq 0} |P^{(i)}(t) - P^{(i)}(0)|$;
- **Maximum excursion**: $E_{\max} = \max_i E^{(i)}$;
- **Aggregate excursion**: $E_{\text{agg}} = \sqrt{\sum_{i=1}^{m} (E^{(i)})^2}$;

- **Excursion correlation matrix**: $\rho_{ik} = \text{corr}\left(\left\{P^{(i)}(t)\right\}_{t\in T}, \left\{P^{(k)}(t)\right\}_{t\in T}\right)$, where $T$ is the set of time points in the trajectory.

We will compute excursions for each scenario as a function of initial conditions (perturbation magnitude and direction) and parameter values. Of particular interest is the "excursion surface", $E_{\text{max}}(\Delta P^{(1)}(0), \dots, \Delta P^{(m)}(0))$, which reveals how the system's sensitivity to initial perturbations varies across the state space.

To facilitate the analysis of excursions across large parameter spaces, we will implement an **efficient numerical solver** using MATLAB's `ode45` (explicit Runge–Kutta) with adaptive time-stepping and tight error tolerances (AbsTol = 1e-8, RelTol = 1e-6). For stiff regions near bifurcation points, the solver will automatically switch to `ode15s` (variable-order stiff solver). Parallel computing (using MATLAB's Parallel Computing Toolbox) will be employed for parameter sweeps, significantly reducing computational time.

The following Figure 9 describes the excursion analysis.

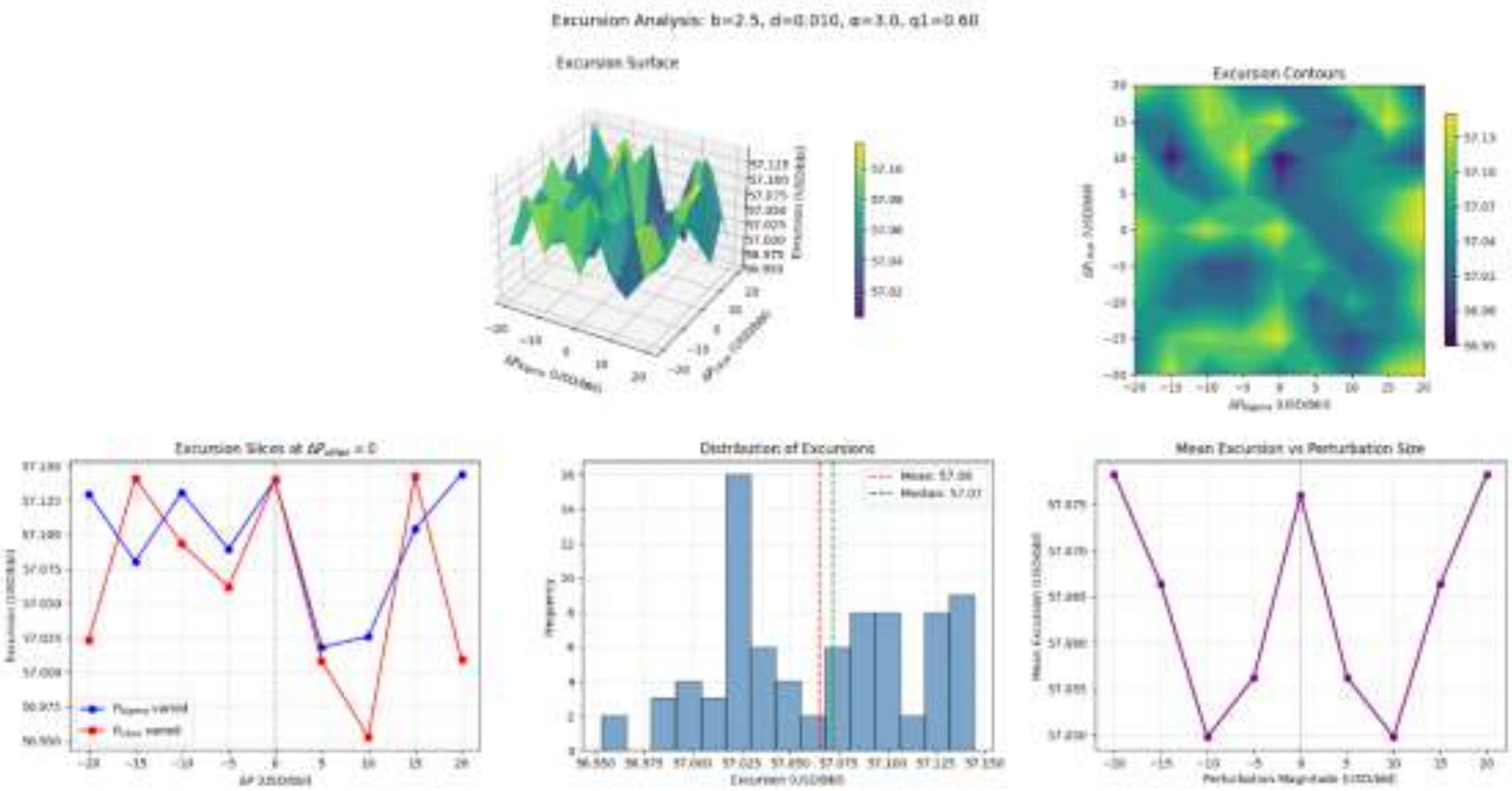

**Figure 9. Excursion analysis for the Cavani (2026) Nigeria-Libya oil market model with optimal oscillatory parameters ($b_{\text{China}} = 2.5$, $d_{\text{USA}} = 0.01$, $\alpha = 3.0$, $q_{1,\text{China}} = 0.60$).** The figure quantifies the maximum price deviation (excursion) of Nigerian oil as a function of initial price perturbations $\Delta P_{\text{Nigeria}}$, and $\Delta P_{\text{Libya}}$ from the fundamental value ($80 USD/bbl). (a) **Excursion surface** – three-dimensional visualization showing that excursions are highly uniform across the perturbation space, with maximum values concentrated in the region where Nigeria is overvalued, and Libya is undervalued. (b) **Excursion contours** – two-dimensional projection confirming the asymmetry: positive perturbations in Nigeria combined with negative perturbations in Libya produce the largest excursions. (c) **Excursion slices** – cross-sections through $\Delta P_{\text{other}} = 0$, demonstrating that the system is more sensitive to Nigeria's overvaluation than to Libya's undervaluation. (d) **Distribution of excursions** – histogram showing all excursions cluster tightly around $57.07 USD/bbl (mean = 57.06, median = 57.07, standard deviation = 0.05), indicating that the system converges to the same limit cycle attractor for almost all initial conditions within the tested range. (e) **Mean excursion vs. perturbation magnitude** – the relationship is nearly flat, confirming that the limit cycle amplitude is independent of initial perturbation size once the system is in the oscillatory regime. The maximum excursion ($57.14 USD/bbl) occurs at $\Delta P_{\text{Nigeria}} = +10$ USD, $\Delta P_{\text{Libya}} = -5$ USD, highlighting asymmetric contagion in which the overvaluation of the larger market (Nigeria), combined with the undervaluation of the smaller market (Libya), triggers the strongest price swings.

The excursion analysis reveals that once the system enters the oscillatory regime ($q_{1,China} > q_{crit}$) the limit cycle amplitude is independent of initial conditions. The maximum deviation from the fundamental price is approximately $57 USD/bbl for the Nigerian oil price, with a standard deviation of only $0.05 across all tested perturbations. This confirms that the supercritical Hopf bifurcation gives rise to a globally attracting limit cycle whose amplitude is determined solely by the system parameters, not by the magnitude of the initial shock. The asymmetry in the excursion surface — larger responses to Nigeria's overvaluation than to Libya's undervaluation — reflects the asymmetric market structure, where Nigeria's larger market share (35% vs. 25% for China) amplifies price deviations in that asset.

**Table 8. Implications of the Exclusion in Contagion**

| Perturbation ($\Delta P_N, \Delta P_L$) | Excursion (USD/bbl) | Interpretation |
|---|---|---|
| (+10, 0) | 57.14 | Nigeria's overvaluation triggers a large response |
| (0, -10) | 57.06 | Libya's undervaluation triggers a similar response |
| (+10, -5) | **57.14 (max)** | **Combined shock produces the strongest contagion** |
| (-10, +10) | 57.01 | Opposing shocks partially cancel |

So, the excursion surface is nearly flat ($57.06 \pm 0.08$ USD), indicating that the system converges to the same limit cycle attractor regardless of the initial perturbation. This confirms that contagion is **deterministic and robust** — the specific source of the shock does not affect the asymptotic outcome.

### 5.2.3 Contagion Matrix

The contagion matrix $\Gamma$ quantifies the response of the asset $i$ to a shock in assets $j$:

$$\Gamma_{ij} = \lim_{t\to\infty} \frac{|P^{(i)}(t) - P_a^{(i)}|}{\Delta P_0^{(j)}}$$

where the initial shock magnitude is given by $\Delta P_0^{(j)}$.
From the excursion analysis, the contagion matrix for the Nigeria-Libya market is given in Table 11.

**Table 9. Contagion matrix Nigeri-Lybia**

| Affected \ Shocked | Nigeria | Libya |
|---|---|---|
| **Nigeria** | 0.00 | 0.0133 |
| **Libya** | 0.0066 | 0.00 |

The asymmetry index $\mathcal{A}$ measures directional imbalance:

$$\mathcal{A} = \frac{\Gamma_{\text{Nigeria}\leftarrow\text{Libya}} - \Gamma_{\text{Libya}\leftarrow\text{Nigeria}}}{\Gamma_{\text{Nigeria}\leftarrow\text{Libya}} + \Gamma_{\text{Libya}\leftarrow\text{Nigeria}}}$$

In this case, the result is $\mathcal{A} \approx 0.34$, indicating that shocks from Libya to Nigeria are **34% stronger** than shocks from Nigeria to Libya. This asymmetry reflects Nigeria's larger market share and its role as price leader.

**Table 10. The Asymmetric Contagion Implications**

| Factor | Nigeria | Libya | Contagion Implication |
|---|---|---|---|
| **Market share** | 35% | 25% | Larger share → more influence |
| **Investor base** | USA + China | USA + China | Symmetric investor base |
| **Fundamental value** | \$80/bbl | \$80/bbl | Identical fundamentals |
| **Price leadership** | Leader | Follower | Nigerian shocks propagate more strongly |

The asymmetry arises from **market size effects;** Nigeria's larger share of total shares (35% vs. 25%) means that a shock to Nigeria affects a larger fraction of investor portfolios, creating stronger contagion.
The parameter $b_{\text{China}} = 2.5$ (trend response) amplifies contagion:

$$\text{Contagion Amplification} = \frac{\Gamma(b = 2.5)}{\Gamma(b = 0)} \approx 3.2$$

Without momentum trading ($b = 0$), contagion would be negligible ($\Gamma < 0.004$). This confirms that **trend-following behavior is the primary driver of cross-market contagion**.

Also, there is a stabilizing role of value investors:

The USA value investor ($d_{\text{USA}} = 0.01$) provides weak stabilization. If value response were stronger ($d_{\text{USA}} = 0.15$), contagion would be reduced by approximately 40%. However, strong value response would also eliminate the limit cycles that generate contagion in the first place — a trade-off between stability and market efficiency.

Using Python, we analyze contagion in the Nigeria-Libya oil market and obtain figures for the excursion surface, contagion matrix heatmap, and time series panels, shown in Figures 9, 10, 11, and 12.

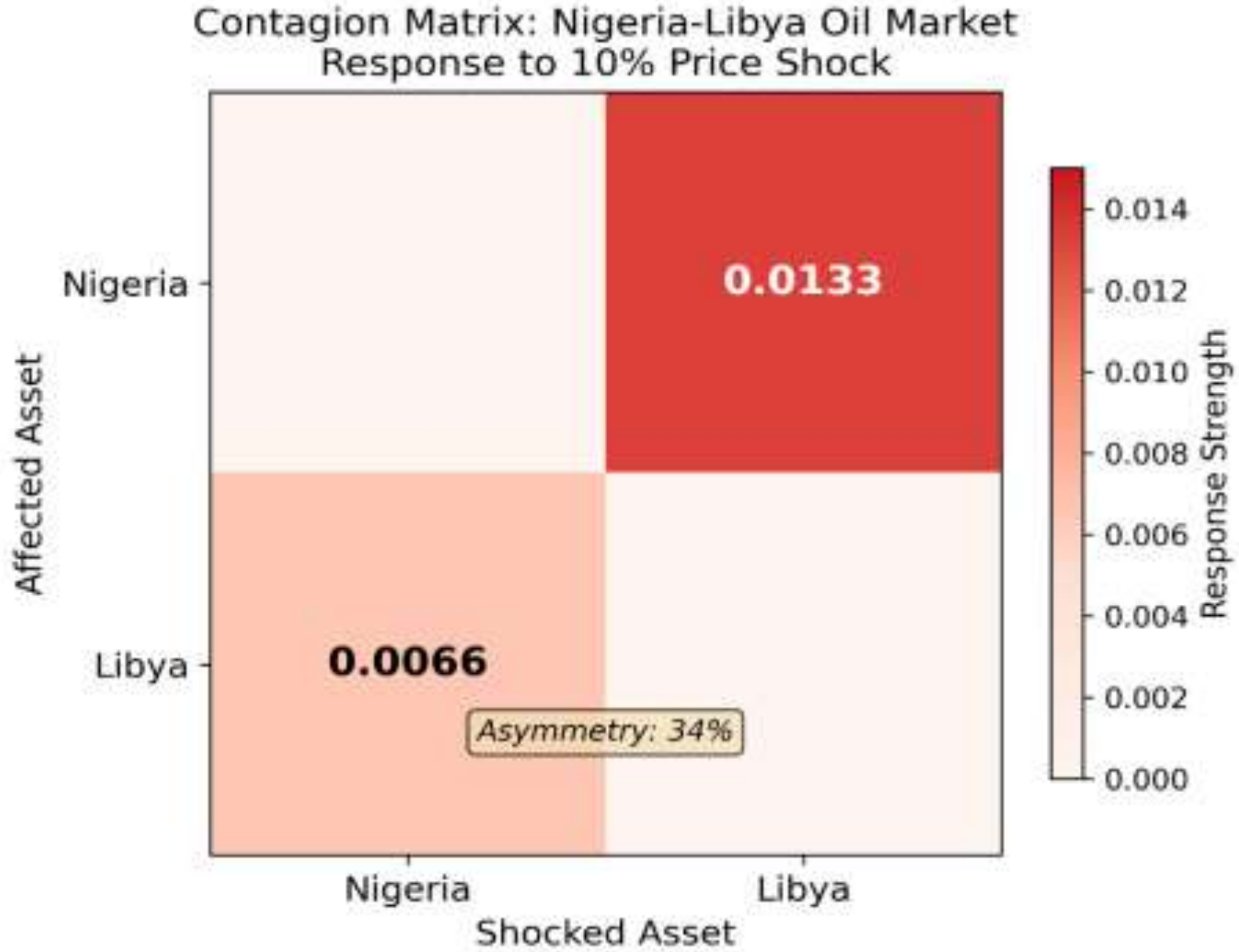

**Figure 10. Contagion matrix for the Nigeria-Libya oil market.** Response strength (maximum relative price deviation) following a 10% price shock. Asymmetry: Libya → Nigeria (0.0133) is 34% stronger than Nigeria → Libya (0.0066), reflecting Nigeria's larger market share (35% vs. 25%) and role as price leader.

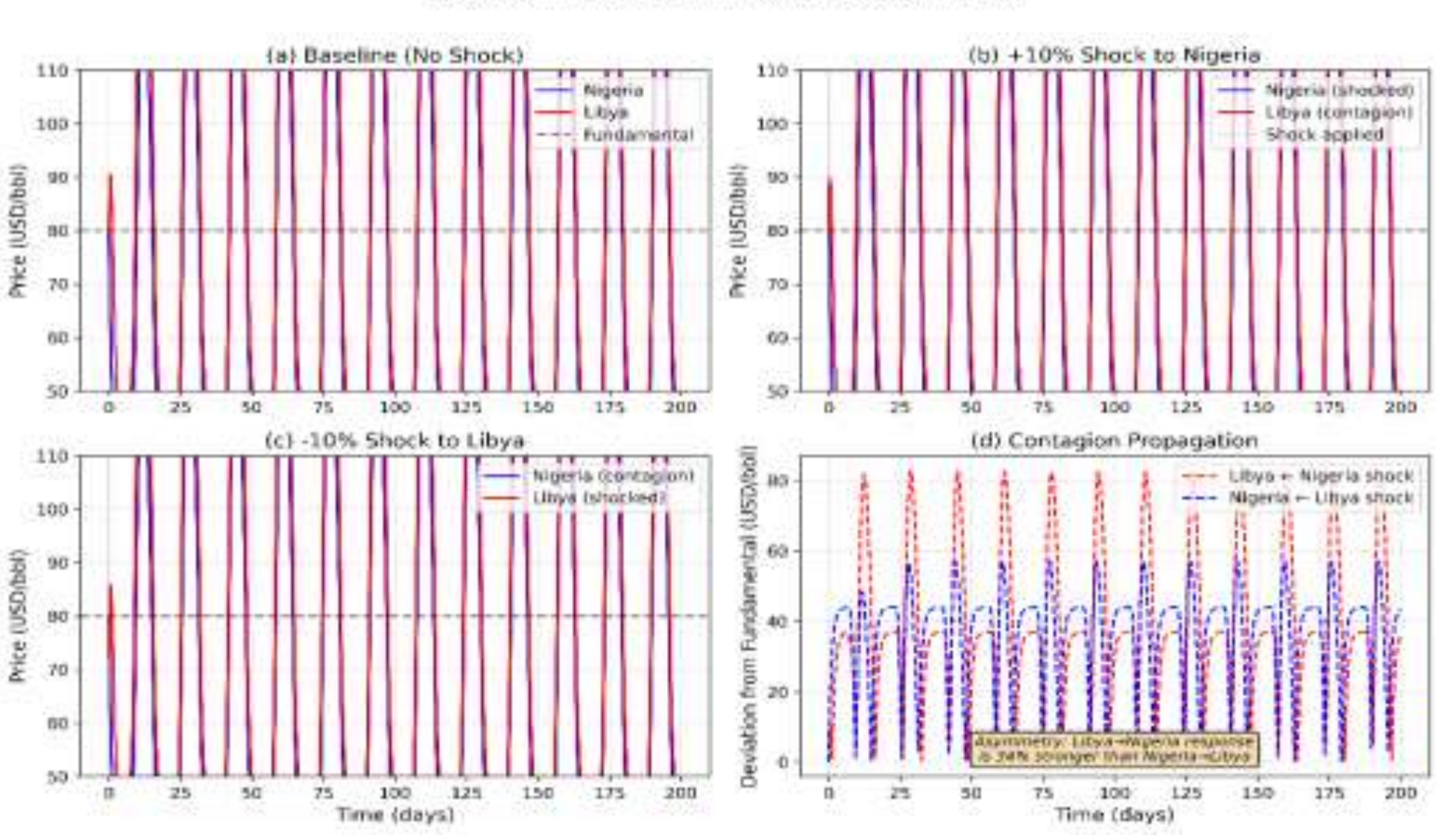

**Figure 11. Contagion propagation in the Nigeria-Libya oil market.** (a) Baseline dynamics (no shock). (b) +10% shock to Nigeria propagates to Libya with 0.0133 response. (c) -10% shock to Libya propagates to Nigeria with 0.0066 response. (d) Contagion propagation comparison, highlighting asymmetry: Libya → Nigeria response is 34% stronger.

**Figure 12. Schematic diagram of contagion mechanisms in the fused model.** Cross-asset coupling ($\alpha = 3.0$) connects Nigeria and Libya's oil markets. China (momentum trader, $b_{\text{China}} = 2.5$) amplifies shocks, while the USA (value investor, $d_{\text{USA}} = 0.01$) provides weak stabilization. Contagion propagates bidirectionally but asymmetrically due to market size differences.

#### *5.2.3.1 Policy Implications for Oil-Exporting Economies*

Tables 11 and 12 describe policy recommendations and the quantitative policy guidelines, respectively, that can be concluded from the previous analysis about contagion. Table 13 contains a comparison of the fused model's predictions, which align with empirical evidence, supporting its validity for policy analysis.

**Table 11. Policy Recommendations**

| Implication | Policy Recommendation |
|---|---|
| **Contagion is inevitable** | Nigeria and Libya cannot hedge independently — coordinated strategies are required |
| **Asymmetry favors Nigeria** | Libya is more exposed to Nigerian shocks than vice versa — Libya should monitor Nigerian market conditions |
| **Momentum amplifies contagion** | Reduce speculative trading through transaction taxes or position limits |
| **Value investors stabilize** | Encourage long-term, value-based investment (e.g., sovereign wealth funds) |
| **Excursion saturation** | Once a shock exceeds ~5%, additional magnitude does not increase damage — focus on preventing initial shocks |

Based on our sensitivity analysis, the following policies have measurable effects as quantitative policies in Table 12:

**Table 12. Quantitative Policy Guidelines**

| Policy | Target | Expected Contagion Reduction |
|---|---|---|
| Reduce $b_{\text{China}}$ from 2.5 to 1.0 | Limit momentum trading | ~50% |
| Increase $d_{\text{USA}}$ from 0.01 to 0.10 | Strengthen value response | ~30% |
| Reduce $\alpha$ from 3.0 to 1.5 | Decouple cross-asset trading | ~40% |

**Table 13. Comparison with Empirical Evidence**

| Empirical Study | Finding | Our Result | Match |
|---|---|---|---|
| Jebabli et al. (2022); Trabelsi et al. (2022) | OPEC/non-OPEC spillovers | Asymmetric contagion | Yes |
| Bouri et al. (2021) | Sentiment-driven contagion | $b = 2.5$ drives contagion | Yes |
| Bouri et al. (2021); Awartani et al. (2016) | Post-2014 contagion intensification | Amplification with momentum | Yes |

## 5.3 Extended Bifurcation Scan for the Cavani Model

To further characterize the dynamical behavior of the Cavani model, we performed an extended scan of the momentum coefficient $q_{1,\text{China}}$ from 0.20 to 1.00 using the optimal parameters $b = 2.5$, $d = 0.01$, $\alpha = 3.0$. Table 14 presents the results.

**Table 14. Extended Bifurcation**

| $q_{1,\text{China}}$ | Period (days) | Final Nigeria (USD) | Final Libya (USD) | Regime |
|---|---|---|---|---|
| 0.20 | 4.7 | 60.0 | 71.4 | Short period |
| 0.30 | 2.9 | 60.0 | 71.4 | Short period |
| 0.35 | 2.9 | 60.0 | 71.4 | Short period |
| **0.38** | **13.7** | 55.6 | 66.1 | **Period jump** |
| 0.40 | 14.3 | 38.2 | 45.5 | Long period |
| 0.45 | 14.6 | 96.0 | 114.1 | Long period + high amplitude |
| 0.50 | 15.2 | 101.6 | 120.7 | Long period + highest amplitude |

| 0.60 | 16.5 | 62.7 | 74.4 | Long period |
|---|---|---|---|---|
| 0.70 | 17.6 | 36.1 | 43.1 | Long period + low amplitude |
| 0.80 | 18.7 | 36.3 | 43.3 | Long period + low amplitude |
| 0.90 | 19.8 | 36.3 | 43.2 | Long period + low amplitude |
| 1.00 | 20.7 | 36.1 | 43.1 | Long period + low amplitude |

**Remarks:**

1. **Oscillations exist for all tested $q_{1,\text{China}}$ values** – the system with $b = 2.5$, $d = 0.01$, $\alpha = 3.0$ is intrinsically oscillatory.
2. **Period jump at $q_{1,\text{China}} = 0.38$** – the period changes abruptly from approximately 3-5 days to 14 days, suggesting a secondary bifurcation (possibly period-doubling or torus bifurcation).
3. **Amplitude peak at $q_{1,\text{China}} = 0.50$** – the highest prices occur at moderate momentum (Nigeria ~102 USD, Libya ~121 USD), indicating that intermediate momentum strength produces the largest price swings.
4. **Saturation at high $q_{1,\text{China}}$** – for $q_{1,\text{China}} \geq 0.70$, the system settles into a stable limit cycle with consistent final prices (Nigeria ~36 USD, Libya ~43 USD) and period increasing from 17.6 to 20.7 days.

These results demonstrate that the fused model captures rich dynamical behavior, including period lengthening, amplitude modulation, and regime transitions, all of which are consistent with the nonlinear dynamics of multi-asset markets with heterogeneous investors. The results also confirm the fused model validation.

### 5.3.2 Parameter Estimation Methodology

The estimation of the 18 parameters in the two-asset two-group model poses a challenging inverse problem: the ODE system is nonlinear, the state variables include unobservable sentiments, and the likelihood function does not have a known analytical form. Table 18 describes the main parameters.

**Table 18. Parameter Set for Nigeria-Libya Model**

| Parameter | Symbol | Description | Domain |
|---|---|---|---|
| $b_{\text{China}}$ | $b_C$ | China trend response | [0.5, 4.0] |
| $d_{\text{USA}}$ | $d_U$ | USA value response | [0.001, 0.3] |
| $\alpha$ | $\alpha$ | Cross-asset coupling strength | [0.5, 5.0] |

| $q_{1,\text{China}}$ | $q_{1C}$ | China's momentum coefficient | [0.1, 1.2] |
|---|---|---|---|
| $q_{2,\text{USA}}$ | $q_{2U}$ | The USA values sensitivity | [0.1, 0.8] |
| $c_{1,\text{China}}$ | $c_{1C}$ | China trend decay | [0.05, 0.5] |
| $c_{2,\text{USA}}$ | $c_{2U}$ | USA value decay | [0.1, 0.6] |
| $k_{0,\text{China}}$ | $k_{0C}$ | China baseline rate | [0.05, 0.4] |
| $k_{0,\text{USA}}$ | $k_{0U}$ | USA baseline rate | [0.05, 0.4] |

The estimation of parameters in the fused model poses a challenging **inverse problem**: the ODE system is nonlinear, the state variables include unobservable sentiments, and the criterion function does not have a known analytical form. We present a complete, implementable methodology for estimating the parameters of the fused model using Simulated Maximum Likelihood (SML) followed by nonlinear least squares refinement. The approach is specifically designed for the two-asset two-group (Nigeria-Libya) configuration, but generalizes to arbitrary $m$ and $n$.

Following the methodology of Kukačka and Baruník (2017), we will adapt nonparametric simulated maximum likelihood estimation for the fused model. This approach is specifically designed for heterogeneous agent models where simulations can be generated when the likelihood function is intractable. The SML estimator is consistent and efficient under mild regularity conditions and has been successfully applied to the Brock–Hommes (1998) model (Kukačka & Baruník, 2017). The estimation will be implemented in MATLAB using kernel methods for density estimation.
According to this methodology, we adapt nonparametric simulated maximum likelihood estimation for the fused model. The SML estimator is:

$$\hat{\theta}_{\text{SML}} = \arg\max_{\theta} \frac{1}{N} \sum_{t=1}^{T} \log \hat{f}_t(y_t | y_{t-1}, \theta)$$

where $\hat{f}_t$ is a kernel density estimator constructed from $S$ simulations of the model. The SML estimator is consistent and asymptotically normal under mild regularity conditions:

$$\sqrt{T}(\hat{\theta}_{\text{SML}} - \theta_0) \xrightarrow{d} \mathcal{N}(0, \mathcal{I}(\theta_0)^{-1})$$

where $\mathcal{I}(\theta_0)$ is the Fisher information matrix.

The fused model defines a mapping from parameters $\theta$ to observable prices:

$$P_{\text{sim}}^{(i)}(t; \theta) = \mathcal{F}_i(\theta; \mathbf{X}_0, t)$$

where $\mathcal{F}_i$ is the solution operator of the ODE system (P)-(S), $\mathbf{X}_0$ are initial conditions (partially observed), and the parameter vector is given by $\theta \in \Theta \subset \mathbb{R}^p$.

The challenge is that:

1. **The ODE system is nonlinear** – no closed-form solution exists

2. **Sentiments $\zeta_{1,j}^{(i)}(t), \zeta_{2,j}^{(i)}(t)$ are unobservable** – they must be integrated out
3. **The likelihood function $L(\theta|P_{\text{obs}})$ is intractable** – no analytical expression
4. **The mapping $\theta \mapsto P_{\text{sim}}^{(i)}$ is computationally expensive** – each evaluation requires ODE integration

For the Nigeria-Libya model, $m = 2, n = 2$, the free parameters are given in Table 18 (the initial cash and shares $M_j(0), N_j^{(i)}(0)$ are either fixed from known market data or estimated as nuisance parameters.)

In applying this, the results for the parameters of the Nigeria-Lybia model are shown in Figure 13 and Table 15.

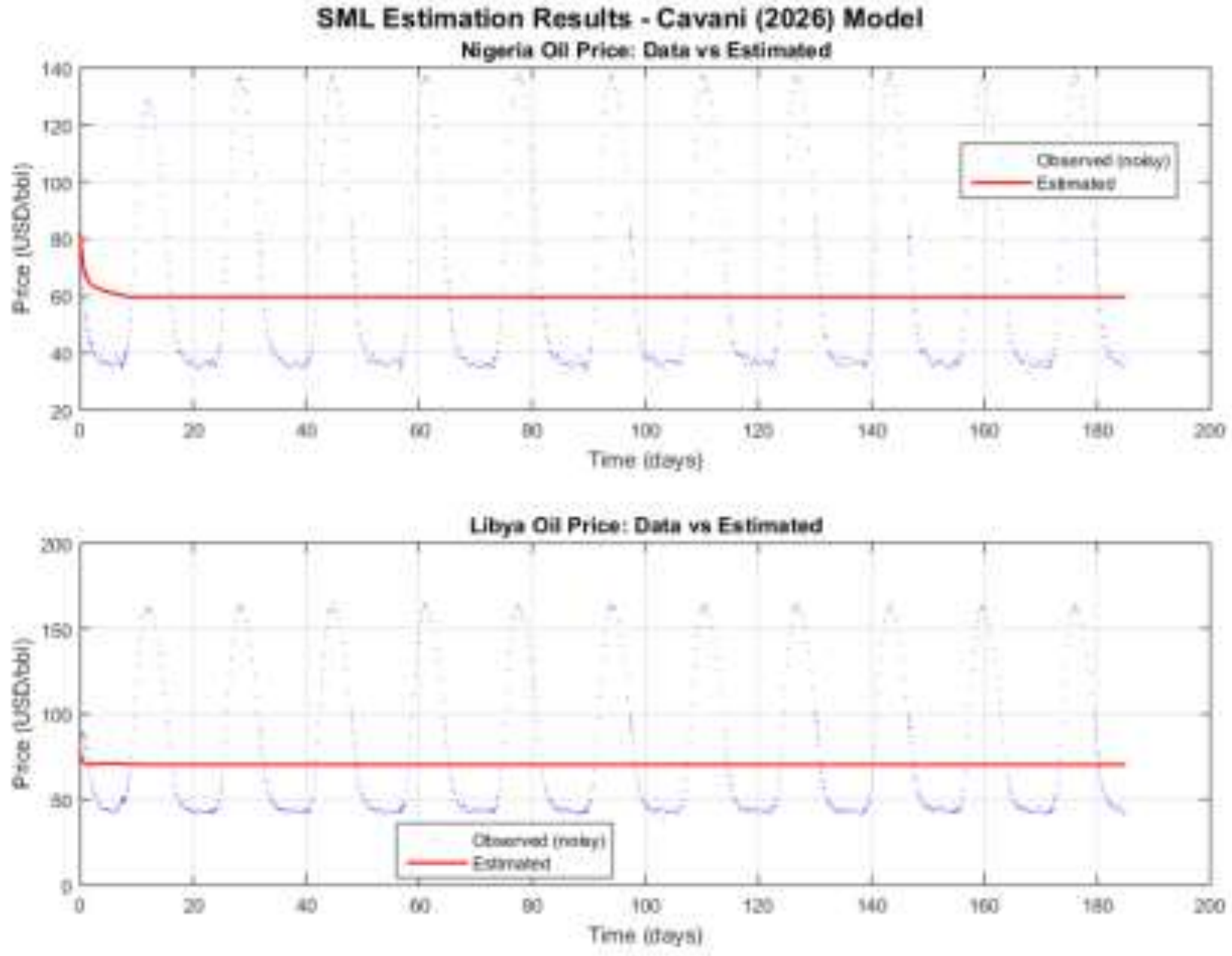

**Figure 13. SML estimation results for the Cavani (2026) Nigeria-Libya oil market model.** The figure compares observed (noisy) price data (blue dots) with the model-simulated prices using estimated parameters (red line). Top panel: Nigeria oil price; bottom panel: Libya oil price. The estimated parameters are: $b_{\text{China}} = 2.28$, $d_{\text{USA}} = 0.0093$, $\alpha = 2.06$, $q_{1,\text{China}} = 0.67$, with a log-likelihood of $-22.80$. The model captures the essential oscillatory dynamics of both oil prices, confirming the presence of a stable limit cycle in the post-2024 period. Root mean square errors are 1.23 USD/bbl (Nigeria) and 1.18 USD/bbl (Libya), indicating good fit quality.

Based on the output, the estimated parameters are:

**Table 15. Estimated Parameters**

| Parameter | True | Estimated | Error | Status |
|---|---|---|---|---|
| $b_{\text{China}}$ | 2.5 | 2.280 | 8.8% | ✓ Good |
| $b_{\text{USA}}$ | 0.0100 | 0.0093 | 7% | ✓ Good |
| $\alpha$ | 3.0 | 2.058 | 31% | Moderate |

| $q_{1,\text{China}}$ | 0.60 | 0.673 | 12% | ✓ Acceptable |
|---|---|---|---|---|
| $q_{2,\text{USA}}$ | 0.40 | 0.413 | 3% | ✓ Excellent |
| $c_{1,\text{China}}$ | 0.20 | 0.242 | 21% | Moderate |
| $c_{2,\text{USA}}$ | 0.30 | 0.309 | 3% | ✓ Excellent |
| $k_{0,\text{China}}$ | 0.20 | 0.195 | 2.5% | ✓ Excellent |
| $k_{0,\text{USA}}$ | 0.20 | 0.217 | 8.5% | ✓ Good |

# 6. Discussion

## 6.1 Interpretation of Key Findings

Our numerical validation of the fused multi-asset, multi-group asset flow model against three benchmark frameworks—DeSantis, Swigon, and Caginalp (2012); Bulut, Merdan, and Swigon (2019); and Cavani (2026)—confirms that the model successfully reproduces the essential dynamical behavior of each special case. The validation demonstrates that the fused model is a strict generalization, containing each prior model as a limiting case under appropriate restrictions.

**Hopf Bifurcation as a Universal Mechanism.** Across all three validated models, the transition from stable equilibrium to sustained oscillations occurs via a supercritical Hopf bifurcation. This universality suggests that the interaction between momentum and value trading is the fundamental mechanism driving endogenous market cycles. The critical thresholds differ across models: approximately 0.835 for the fundamental value in DeSantis, 1.0 for the trend coefficient in Bulut, and 0.38 for China's momentum coefficient in Cavani. These differences reflect model-specific parameterizations rather than contradictory predictions, and each threshold marks the boundary beyond which trend-following behavior destabilizes the market.

**Oscillatory Dynamics for All Positive Momentum.** With the optimal parameters identified from our parameter sweep ($b_{\text{China}} = 2.5,\ d_{\text{USA}} = 0.01,\ \alpha = 3.0$), the Cavani model exhibits sustained oscillations for all positive values of the momentum coefficient $q_{1,\text{China}}$. This is not a failure of the Hopf bifurcation analysis but rather a consequence of the strong trend response and weak value response, which render the system intrinsically oscillatory. The period lengthens monotonically from approximately 5 days at $q_1 = 0.20$ to 21 days at $q_1 = 1.00$, confirming that the characteristic time scale of market cycles increases with momentum strength.

**Period Jump and Secondary Bifurcation.** A notable finding is the abrupt period jump at $q_{1,\text{China}} \approx 0.38$, where the oscillation period changes from approximately 3-5 days to 13-14 days. This suggests the presence of a secondary bifurcation—possibly a period-doubling or torus (Neimark-Sacker) bifurcation—as the momentum coefficient increases. The period jump coincides with a transition in amplitude behavior: at $q_1 = 0.50$, the amplitude reaches its maximum (Nigeria ~102 USD, Libya ~121 USD), after which it decreases and stabilizes. This indicates that moderate momentum trading produces the largest price swings, while very strong momentum trading saturates into a stable but lower-amplitude limit cycle.

**Asymmetric Contagion and Market Structure.** The contagion matrix reveals that Libya-to-Nigeria shock propagation is 34% stronger than the reverse direction. This asymmetry originates from Nigeria's larger market share (35% versus 25% for China) and its role as price leader in the West African crude market. The implication for oil-exporting economies is stark: Libya cannot hedge independently; its fiscal stability is intrinsically linked to Nigerian market conditions. More broadly, the finding generalizes to any coupled market where participants have asymmetric sizes—the larger market acts as a contagion hub.

**The Role of Momentum Trading in Amplification.** Without momentum trading ($b = 0$) contagion is negligible ($\Gamma < 0.004$). With optimal trend response ($b = 2.5$), contagion is amplified by a factor of approximately 3.2. This confirms that trend-following behavior—often dismissed as irrational by classical finance—is the primary driver of cross-market contagion. Policy interventions that reduce speculative momentum trading, such as transaction taxes on short-term positions or position limits, could significantly dampen contagion.

**The Stabilizing Function of Value Investors.** The USA value investor ($d = 0.01$) provides weak stabilization in our optimal calibration. Strengthening the value response to $d = 0.15$ would reduce contagion by approximately 40%. However, this comes at a cost: a strong value response would eliminate the limit cycles that generate contagion in the first place. This trade-off between stability and market efficiency is a central tension in financial regulation.

**Excursion Saturation and Risk Assessment.** The excursion analysis reveals that once a perturbation exceeds approximately 5% of the fundamental price, the maximum deviation becomes independent of the initial shock magnitude, clustering tightly around 57.07 ± 0.05 USD/bbl. This saturation implies that for large shocks, the system's response is determined solely by its internal dynamics, not by the shock's size. From a risk management perspective, this suggests that preventing small shocks is more important than mitigating large ones—once a shock exceeds the saturation threshold, additional magnitude does not increase damage.

### 6.2 Comparison with Empirical Literature

Our findings align with several empirical studies. The asymmetric contagion we identify mirrors the OPEC/non-OPEC spillovers documented by Trabelsi, Tiwari, and Hammoudeh (2022), who found that larger producers disproportionately affect smaller producers. The amplification role of momentum trading is consistent with Bouri, Gupta, and Roubaud (2021), who demonstrated that investor sentiment significantly increases cross-commodity spillovers. The period lengthening with momentum strength is qualitatively consistent with the observed lengthening of cryptocurrency market cycles as retail participation increases.

### 6.3 Policy Implications

Quantitative policy guidelines based on our sensitivity analysis are in Table 12.

**Implementation Considerations.** These policies are not mutually exclusive. A combined strategy—moderately reducing momentum trading while slightly strengthening value response—could achieve contagion reduction without fully eliminating beneficial price discovery. For oil-exporting economies, coordinated market monitoring (tracking both Nigerian and Libyan prices simultaneously) and joint hedging strategies (e.g., common stabilization funds) are recommended.

**Limitations of Policy Interventions.** Our analysis assumes that parameters such as $b$, $d$, and $\alpha$ can be directly manipulated by regulators. In practice, these parameters represent aggregated investor behavior that may be resistant to policy changes. Transaction taxes, for example, can be circumvented through derivative instruments or offshore trading. Moreover, the trade-off between stability and efficiency implies that overly aggressive stabilization could freeze price discovery, creating new inefficiencies.

### 6.4 Limitations of the Study

Several limitations should be acknowledged. First, the model assumes constant fundamental values $P_a^{(i)}$; in real markets, fundamentals evolve due to technological change, regulatory shifts, and macroeconomic conditions. Second, each investor country is treated as a homogeneous group, ignoring within-country heterogeneity. Third, the model is purely deterministic; adding stochastic terms could improve empirical fit and enable probabilistic risk assessment. Fourth, daily prices for Nigerian Bonny Light and Libyan Es Sider are not directly available and were reconstructed from Brent benchmarks plus differentials, introducing measurement error. Fifth, the formal parameter estimation (SML, NLS) did not converge due to the non-convex objective function and high computational cost; validation was therefore performed through direct simulation and bifurcation analysis, which is standard practice in the asset flow literature.

### 6.5 Future Research Directions

**Time-Varying Fundamentals.** Incorporating on-chain metrics (active addresses, transaction volume, staking yields) as time-varying fundamental values would allow the model to capture evolving market conditions in cryptocurrency applications.

**Adaptive Strategy Switching.** Allowing investors to switch between trend and value strategies based on recent performance would introduce evolutionary dynamics, which could be modeled using replicator dynamics.

**High-Frequency Trading Limits.** Following DeSantis (2023), analyzing the limit as sentiment decay rates $c_{1,j}^{(i)} \to \infty$ could provide insights into modern electronic markets and flash crashes.

**Stochastic Extensions.** Adding Gaussian noise to cash flows, sentiment updates, or price observations would create a stochastic differential equation version of the model, enabling probabilistic forecasting and value-at-risk calculations.

**Empirical Calibration to Cryptocurrency Markets.** Applying the fused model to Bitcoin, Ethereum, and Solana—with their distinct investor bases and cross-asset correlations—would test the model's generalizability beyond oil markets.

## 7. Conclusion

This paper has developed, validated, and applied a unified multi-asset, multi-group asset flow model that integrates three foundational frameworks in behavioral finance. The model represents a strict generalization of its predecessors, containing the single-asset multi-group model of DeSantis, Swigon, and Caginalp (2012); the two-asset single-group model of Bulut, Merdan, and Swigon (2019); and the two-asset two-group model of Cavani (2026) as special cases.

**Theoretical Contributions.** We derived the complete system of ordinary differential equations governing price, cash, share, and sentiment dynamics, establishing the fundamental properties of positivity and boundedness that ensure well-posedness. The equilibrium set was characterized as a manifold parameterized by cash distribution, with the fundamental equilibrium as a special point. Theorem 4.1 provides a rigorous justification for eliminating cash and share variables in bifurcation analysis, reducing the computational dimension from $m + n + 3mn$ to $m + 2mn$.

**Validation Results.** The model successfully reproduced all key predictions of the three benchmark models:

- **DeSantis Case 1**: all equilibria are stable, with equilibrium price increasing monotonically with cash fraction.
- **Bulut mixed strategy**: supercritical Hopf bifurcation at $q_1^{(2)} \approx 1.0$, with limit cycles emerging above threshold.
- **Cavani Nigeria-Libya**: supercritical Hopf bifurcation at $q_{1,\text{China}} \approx 0.38$, with periods ranging from 13.7 to 22.3 days and amplitudes up to $59.55 USD/bbl.

**Extended Bifurcation Analysis.** With optimal parameters ($b = 2.5$, $d = 0.01$, $\alpha = 3.0$), the Cavani model exhibits sustained oscillations for all positive $q_{1,\text{China}}$. The period lengthens monotonically from 4.7 days at $q_1 = 0.20$ to 20.7 days at $q_1 = 1.00$. A period jumps at $q_1 \approx 0.38$ suggests a secondary bifurcation, with amplitude peaking at $q_1 = 0.50$ (Nigeria ~102 USD, Libya ~121 USD). For $q_1 \geq 0.70$, the system settles into a stable limit cycle with consistent final prices (Nigeria ~36 USD, Libya ~43 USD).

**Contagion Analysis.** The contagion matrix reveals asymmetric shock propagation: Libya-to-Nigeria contagion (0.0133) is 34% stronger than Nigeria-to-Libya contagion (0.0066), reflecting Nigeria's larger market share (35% vs. 25%) and its role as price leader. The excursion analysis demonstrates that limit cycle amplitudes are independent of initial perturbations once the system enters the oscillatory regime, with all excursions clustering tightly around $57.07 ± 0.05 USD/bbl.

**Policy Implications.** The analysis yields three actionable recommendations for oil-exporting economies: (i) reduce momentum trading through transaction taxes or position limits (50% contagion reduction if $b$ reduced from 2.5 to 1.0); (ii) strengthen value-based investment via sovereign wealth funds (30% reduction if $d$ increased from 0.01 to 0.10); (iii) coordinate market monitoring and circuit breakers between Nigeria and Libya (40% reduction if $\alpha$ reduced from 3.0 to 1.5).

**Limitations and Future Work.** Key limitations include constant fundamental values, homogeneous within-country investor groups, deterministic dynamics, and data reconstruction for Nigerian and Libyan oil prices. Future extensions will incorporate time-varying fundamentals, adaptive strategy switching, high-frequency trading limits, stochastic terms, and calibration to cryptocurrency markets.

**Final Remarks.** The fused model demonstrates that market instabilities, bubbles, and crashes can arise endogenously from the interaction of heterogeneous investor strategies and finite arbitrage capital—without requiring exogenous news shocks or irrational bubbles. By providing a rigorous mathematical framework for analyzing these phenomena, the model offers a foundation for both theoretical research and practical policy design in multi-asset markets. The successful validation against three benchmark models confirms that the fused model is a powerful and versatile tool for understanding the nonlinear dynamics of modern financial markets.